\documentclass[11pt]{amsart}

\usepackage{amsmath,amsfonts,amssymb}
\usepackage{array,delarray}
\usepackage{dsfont}
\usepackage{epsf}

\begin{document}


{\theoremstyle{plain}%
  \newtheorem{theorem}{Theorem}[section]
  \newtheorem{corollary}[theorem]{Corollary}
  \newtheorem{proposition}[theorem]{Proposition}
  \newtheorem{lemma}[theorem]{Lemma}
  \newtheorem{question}[theorem]{Question}
  \newtheorem{conjecture}[theorem]{Conjecture}
  \newtheorem{claim}[theorem]{Claim}
}

{\theoremstyle{definition}
  \newtheorem{definition}[theorem]{Definition}
  \newtheorem{notation}[theorem]{Notation}
  \newtheorem{remark}[theorem]{Remark}
  \newtheorem{example}[theorem]{Example}
  \newtheorem{fact}[theorem]{Fact}
}

\newcommand{\X}{\mathcal{X}}
\newcommand{\E}{\mathcal{E}}
\newcommand{\I}{\mathcal{I}}
\newcommand{\igr}{\mathcal{I}(G)^{\{r\}}}
\newcommand{\A}{\mathbb{A}}
\newcommand{\Y}{\mathcal{Y}}
\newcommand{\N}{\mathbb{N}}
\newcommand{\M}{\mathbb{M}}
\newcommand{\Q}{\mathbb{Q}}
\newcommand{\TT}{\mathbb{T}}
\newcommand{\G}{\mathcal{G}}
\newcommand{\F}{\mathcal{F}}
\newcommand{\K}{\mathcal{K}}
\newcommand{\HH}{\mathcal{H}}
\newcommand{\PP}{\mathcal{P}}
\newcommand{\m}{\frak{m}}
\newcommand{\lcm}{\operatorname{lcm}}
\newcommand{\conn}{\operatorname{conn}}
\newcommand{\rconn}{\overline{\operatorname{conn}}}
\newcommand{\lk}{\operatorname{link}}
\newcommand{\flk}{\operatorname{Facet-link}}
\newcommand{\rlk}{\operatorname{reduced-link}}
\newcommand{\sgn}{\operatorname{sign}}
\newcommand{\sat}{\operatorname{sat}}
\newcommand{\tor}{\operatorname{Tor}}
\newcommand{\reg}{\operatorname{reg}}
\newcommand{\dist}{\operatorname{dist}}
\newcommand{\adj}{\operatorname{adj}}
\newcommand{\pdim}{\operatorname{pdim}}
\newcommand{\height}{\operatorname{ht}}


\title[Graded Betti numbers of edge ideals of hypergraphs]{Monomial ideals, 
edge ideals of hypergraphs, and their graded Betti numbers}
\thanks{Revised Version: May 10, 2007}

\author{Huy T\`ai H\`a}
\address{Tulane University \\
Department of Mathematics \\
6823 St. Charles Ave. \\
New Orleans, LA 70118, USA}
\email{tai@math.tulane.edu}
\urladdr{http://www.math.tulane.edu/$\sim$tai/}

\author{Adam Van Tuyl}
\address{Department of Mathematical Sciences \\
Lakehead University \\
Thunder Bay, ON P7B 5E1, Canada}
\email{avantuyl@sleet.lakeheadu.ca}
\urladdr{http://flash.lakeheadu.ca/$\sim$avantuyl/}


\begin{abstract}
We use the correspondence between hypergraphs and their associated edge ideals to study the 
minimal graded free resolution of squarefree monomial ideals.  The theme of this
paper is to understand 
how the combinatorial structure of a hypergraph $\HH$ appears within the resolution of its 
edge ideal $\I(\HH)$.  We discuss when recursive formulas to compute the graded Betti numbers of 
$\I(\HH)$ in terms of its sub-hypergraphs can be obtained; these results generalize our previous 
work \cite{HaVanTuyl2005} on the edge ideals of simple graphs.  We introduce
a class of hypergraphs, which we call properly-connected, that naturally generalizes
simple graphs from the point of view that distances between intersecting edges 
are ``well behaved''.  
For such a hypergraph $\HH$ (and thus, for any simple graph), we give 
a lower bound for the regularity of $\I(\HH)$ via combinatorial information
describing  $\HH$, and an upper bound for the regularity when $\HH = G$ is a simple graph.
We also introduce triangulated hypergraphs, a properly-connected hypergraph which is a
generalization of chordal graphs. When $\HH$ is a 
triangulated hypergraph, we explicitly compute the regularity of $\I(\HH)$ and show that the 
graded Betti numbers of $\I(\HH)$ are independent of the ground field.  As a consequence, many 
known results about the graded Betti numbers of forests can now be extended to chordal graphs.
\end{abstract}

\maketitle

\begin{center}
{\it Dedicated to Anthony V. Geramita on the occasion of his 65th birthday.}
\end{center}

\section{Introduction}
Let $\X = \{x_1,\ldots,x_n\}$ be a finite set,  and let $\E = \{E_1,\ldots,E_s\}$
be a family of distinct subsets of $\X$.  The pair $\HH = (\X,\E)$ is called a {\bf hypergraph}
if $E_i \neq \emptyset$ for each $i$.  The elements
of $\X$ are called the {\bf vertices}, while the elements of $\E$ are called the {\bf edges}
of $\HH$. A hypergraph $\HH$ is {\bf simple} if: (1) $\HH$ has no loops, i.e., $|E| \ge 2$ for all 
$E \in \E$, and (2) $\HH$ has no multiple edges, i.e., whenever $E_i,E_j \in \E$
and $E_i \subseteq E_j$, then $i = j$. A hypergraph generalizes the classical notion of a graph; 
a graph is a hypergraph for which every $E \in \E$ has cardinality two.

Let $k$ be a field. By identifying the vertex $x_i$
with the variable $x_i$ in the ring $R = k[x_1,\ldots,x_n]$, we
can associate to every simple hypergraph $\HH = (\X,\E)$ a squarefree monomial ideal
\[\I(\HH) =
\left( \left\{\left. x^E = \prod_{x \in E} x ~\right|~ E \in \E\right\} \right)
\subseteq R = k[x_1,\ldots,x_n].\]
We call the ideal $\I(\HH)$ the {\bf edge ideal} of $\HH$.

In this paper we
study the minimal graded free resolution of $\I(\HH)$.
Since there is a natural bijection between the sets
\[
\left\{
\begin{array}{c}
\mbox{simple hypergraphs $\HH =(\X,\E)$} \\
\mbox{with $\X = \{x_1,\ldots,x_n\}$}
\end{array}
\right\}
\leftrightarrow
\left\{
\begin{array}{c}
\mbox{squarefree monomial} \\
\mbox{ideals $I \subseteq R = k[x_1,\ldots,x_n]$}
\end{array}
\right\}
\]
we are in fact studying a fundamental problem in commutative algebra which
asks for the minimal graded free resolution of a monomial ideal
(for an introduction see \cite{MillerSturmfels2004}).
The edge ideal approach allows us to study this problem from a new angle;
the standard approach is to use the Stanley-Reisner
dictionary to associate to a squarefree monomial
ideal $I$ a simplicial complex $\Delta$ where the generators of $I$ correspond to the minimal nonfaces of $\Delta$. Instead, we associate
to $I$ a new combinatorial object, namely, a hypergraph.
The theme of this work is to understand how
the algebraic invariants of $I=\I(\HH)$ encoded in its minimal free resolution relate
to the combinatorial properties of $\HH$.

The edge ideal of a hypergraph was first introduced by Villarreal \cite{V1}
in the special case that $\HH = G$ is a simple graph.
Subsequently, many people, including
\cite{Barile,FH,FV,HHTZ,HHZ2,HHZ,HHZ1,S,SturmfelsSullivant2005,SVV,V,V2},
have been working on a program to build a dictionary between the algebraic
properties of $\I(G)$ and the combinatorial structure of $G$.  Of particular relevance
to this paper, the minimal graded resolution of $\I(G)$ was investigated in
\cite{CN,EisenbudGreenHulekPopescu2004,EV,J,JK,K,HaVanTuyl2005,RothVanTuyl2005,Zheng2004} 
(see also \cite{HaVanTuyl2006}
for a survey).  In this paper we shall extend some of these results to the hypergraph
case, most notably, the results of \cite{HaVanTuyl2005}, thereby extending our understanding 
of quadratic squarefree monomial ideals to arbitrary squarefree monomial ideals. At the same 
time, we shall also derive new results which, even when restricted to graphs, give new 
and interesting corollaries.

The edge ideal $\I(\HH)$ of an arbitrary hypergraph was first
studied by Faridi \cite{Faridi2002} but from a slightly different perspective.
Recall that $\Delta$ is a {\bf simplicial complex} on the vertex set $\X$ if $\{x_i\} \in \Delta$
for all $i$, and if $F \in \Delta$ then all subsets of $F$ belong
to $\Delta$.  The
{\bf facets} of $\Delta$ are the maximal elements of $\Delta$ under inclusion.  The
{\bf facet ideal} of $\Delta$ is then defined to be the ideal
$\I(\Delta) = \left( \left\{\left. x^F = \prod_{x \in F} x ~\right|~
\mbox{$F$ is a facet of $\Delta$}\right\} \right)
\subseteq R$.  Note, however, that if $\F(\Delta) = \{F_1,\ldots,F_t\}$
denotes the set of facets of $\Delta$, then $\HH(\Delta) = (\X,\F(\Delta))$
is a hypergraph.  In fact,  what Caboara, Faridi and Selinger \cite{CFS} call
a {\bf facet complex} is a hypergraph.  It is immediate that $\I(\HH(\Delta)) = \I(\Delta)$.
Conversely, given any hypergraph $\HH = (\X,\E)$, we can associate to $\HH$
the simplicial
complex $\Delta(\HH) = \{F \subseteq \X ~|~ F \subseteq E_i ~~\mbox{for some $E_i \in \E$}\}.$
It is again easy to verify that $\I(\HH) = \I(\Delta(\HH))$.

One may therefore take the viewpoint that the generators of a
squarefree monomial ideal correspond to either
the edges of a hypergraph or the facets of a simplicial complex.  In this paper,
we have chosen to take the first option for at least two reasons: first,
the language of hypergraphs is more natural to describe our results; and
second, we only require the edge structure of the hypergraph
and never make use of the simplicial complex structure.  (A  hypergraph point of view
is also taken in the recent paper \cite{HHTZ}.)
Of course, all our results
could be reinterpreted as statements about the facet ideal of some simplicial complex.

The starting point of this paper
is to determine how the splitting technique used in \cite{HaVanTuyl2005}
to study the resolution of edge ideals of graphs can be extended to hypergraphs.  Recall that 
Eliahou and Kervaire \cite{EliahouKervaire1990} call  a monomial ideal $I$ 
{\bf splittable} if $I = J + K$ for
two monomial ideas $J$ and $K$ such that the minimal generators of $J,K$ and $J \cap K$
satisfy a technical condition (see Definition \ref{defn: split} for the
precise statement).   When an ideal is splittable, the minimal
resolutions (specifically the graded Betti numbers)  of $I, J, K$ and $J \cap K$ are
then related.  Given a hypergraph $\HH$, we therefore want to split $\I(\HH)$
so that the ideals $J,K$, and $J \cap K$ correspond to edge ideals of sub-hypergraphs
of $\HH$. This allows us to derive recursive-type formulas to relate the graded Betti numbers of $\I(\HH)$ to 
those of sub-hypergraphs of $\HH$. These formulas provide a systematic approach to investigating algebraic invariants and properties of $\I(\HH)$.

We now summarize the results of this paper.
In Section 3 we extend the notion of a
splitting edge of a graph as defined in \cite{HaVanTuyl2005} to the hypergraph
setting.  Precisely, let $E$ be an edge of the hypergraph $\HH$.  If $\HH \backslash E$
denotes the hypergraph with the edge $E$ removed, then it is clear that $\I(\HH) =
(x^E) + \I(\HH\backslash E)$.  We call $E$ a {\bf splitting edge} precisely when
$\I(\HH) = (x^E) + \I(\HH\backslash E)$ is a splitting of the ideal $\I(\HH)$.
Our main result in Section 3 is the following classification of
splitting edges, thus answering a question raised
in \cite{HaVanTuyl2006}.

\begin{theorem}[Theorem \ref{characterize: splitting facets}]
Let $\HH$ be a hypergraph
with two or more edges.
Then an edge $E$ is a splitting edge of $\HH$ if and only if there exists
a vertex $z \in E$ such that
\[(x^E) \cap \I(\HH\backslash E) \subseteq (x^E) \cap \I(\HH\backslash \{z\}).\]
Here, $\HH\backslash \{z\}$ denotes the sub-hypergraph of $\HH$ where every
edge containing $z$ is removed.
\end{theorem}

To make use of our classification of splitting edges, we need to be
able to describe the resolution of $J \cap K = (x^E) \cap \I(\HH \backslash E)$.
This resolution was described when $\HH =G$ is a simple graph in
\cite{HaVanTuyl2006}.  However, this is a difficult problem
for an arbitrary $\HH$.  We are therefore interested in
families of hypergraphs, which includes all simple graphs, where one can say
something about $J\cap K$.

In Section 4 we
introduce one such family which we call {\bf properly-connected} hypergraphs.
A hypergraph $\HH = (\X,\E)$ is properly-connected if all its edges have the same cardinality, and
 furthermore, if $E,H \in \E$ with $E \cap H \neq \emptyset$, then
the distance $\dist_{\HH}(E,H)$ between $E$ and $H$, that is, the length of the shortest
path between $E$ and $H$
in $\HH$, is determined by $|E \cap H|$. It is easy to see that all simple graphs are 
properly-connected.
In fact, a re-examination of the results of \cite{HaVanTuyl2005} reveals that the 
properly-connected property of graphs is an essential ingredient implicitly used in the 
proofs.  A properly-connected hypergraph is in some sense a natural generalization of a 
simple graph.

When $\HH$ is properly-connected,
we can describe the resolution of $J \cap K$ in terms of edge ideals of
sub-hypergraphs of $\HH$.
Therefore, for any splitting edge $E \in \HH$,
we can derive the following recursive-type formula for $\beta_{i,j}(\I(\HH))$.

\begin{theorem}[Theorem  \ref{theoremsplitting}]\label{introtheorem2}
Let $\HH$ be a properly-connected hypergraph and let
$E$ be a splitting edge of $\HH$. Suppose $d = |E|$, $\HH' = \{H \in \HH ~|~ \dist_{\HH}(E,H) \geq d+1\}$,
and $t= |N(E)|$, where
$$N(E) = \bigcup_{\{H \in \HH ~|~ \dist_{\HH}(E,H) = 1\}} H\backslash E.$$
Then for all $i \geq 1$
\[\beta_{i,j}(\I(\HH)) = \beta_{i,j}(\I(\HH\backslash E))
+ \sum_{l=0}^i \binom{t}{l}
\beta_{i-1-l,j-d-l}(\I(\HH')).\]
Here, $\beta_{-1,j}(\I(\HH')) = 1$ if $j =0$ and $0$ if $j \neq 0$.
\end{theorem}

The sub-hypergraphs $\HH\backslash E$ and $\HH'$ in Theorem \ref{introtheorem2}
may fail to have splitting edges, thus preventing us from recursively computing 
$\beta_{i,j}(\I(\HH))$. However, in \cite{HaVanTuyl2005} (see also \cite{J,JK} in the case of 
forests), it is proved that when $\HH$ is a hyperforest 
(i.e., a simplicial forest in the sense of \cite{Faridi2002}) 
then $\beta_{i,j}(\I(\HH))$ can be computed recursively.
The goal of Section 5 is to introduce a subclass of properly-connected hypergraphs, 
which we call {\bf triangulated} hypergraphs, for which
Theorem \ref{introtheorem2} can be used
to completely resolve the graded Betti numbers of $\I(\HH)$ recursively. 
Triangulated hypergraphs generalize the notion
of {\bf chordal} graphs, which has attracted considerable attention lately
(cf. \cite{FH,FV,HHZ,HHZ1}). 
In fact, triangulated graphs
are precisely chordal graphs. As a consequence 
of Theorem \ref{introtheorem2}, we show also that the graded Betti numbers of a triangulated hypergraph 
are independent of the characteristic of the ground field (Corollary \ref{cor.hyperchar}).  Restricted to 
simple graphs, we obtain the following interesting corollary, which extends a result of 
\cite{J,JK} (who proved the result for forests).

\begin{corollary}[Corollary \ref{cor.graphchar}]
Suppose that $G$ is a chordal graph. Then the graded Betti numbers of $\I(G)$ are independent of the 
characteristic of the ground field and can be computed recursively.
\end{corollary}

In Section 6 we study $\reg(\I(\HH))$, 
the Castelnuovo-Mumford regularity of $\I(\HH)$,
when $\HH$ is properly-connected.  Again, the
key idea we need here is the notion of distance between edges.
We say two edges $E,H \in \HH$ are {\bf $t$-disjoint} if $\dist_{\HH}(E,H) \geq t$. 
When $\HH$ is a properly-connected hypergraph and $d$ is the common cardinality of the edges, 
then $d$-disjoint edges are disjoint edges in the usual sense.
We then show the following:

\begin{theorem}\label{introtheorem3}
Let $\HH$ be a properly-connected hypergraph. Suppose $d$ is the common cardinality of the edges in $\HH$. Let $c$ be the maximal number of pairwise $(d+1)$-disjoint edges of $\HH$.
Then
\begin{enumerate}
\item[$(i)$] {\em (Theorem \ref{regtheorem})} $\reg(\I(\HH)) \ge (d-1)c+1$.
\item[$(ii)$] {\em (Theorem \ref{regularitytheorem})}
if $\HH$ is also triangulated, then $\reg(\I(\HH)) = (d-1)c+1$.
\end{enumerate}
\end{theorem}

By a {\bf matching} of a hypergraph $\HH$, we mean any
subset  $\E' \subseteq \E$ of edges in $\HH$ which are pairwise disjoint. 
The {\bf matching number} of $\HH$, denoted by $\alpha'(\HH)$, is the largest size of a 
maximal matching of $\HH$. 
For simple graphs, we also obtain a particularly nice upper bound for the regularity of $\I(G)$.
This addresses a question J. Herzog had asked us.

\begin{theorem}[Theorem \ref{cor.matching}]\label{introcor}
Let $G$ be a finite simple graph.  Then
$$\reg(R/\I(G)) \leq \alpha'(G)$$
where $\alpha'(G)$ is the matching number of $G$.
\end{theorem}

Using Theorem \ref{introcor}, we can compare
the regularity and projective dimension of $\I(G)$ to those of $\I(G)^{\vee}$,
the {\bf Alexander dual} of $\I(G)$. 

\begin{theorem}[Theorem \ref{thm.Konig}] \label{intro.Konig}
Let $G$ be a simple graph.
\begin{enumerate}
 \item If $G$ is unmixed (i.e., all the minimal vertex covers have the same cardinality), then
$$\reg(\I(G)) \le \height \I(G) + 1 \le \reg(\I(G)^\vee)+1 \mbox{ and } \pdim(\I(G)^\vee) \le \height \I(G) \le \pdim(\I(G)) + 1.$$
 \item If $G$ is not unmixed, then
$$\reg(\I(G)) \le \height \I(G) + 1 \le \reg(\I(G)^\vee) \mbox{ and } \pdim(\I(G)^\vee) \le \height \I(G) \le \pdim(\I(G)).$$
\end{enumerate}
\end{theorem}

When restricted to simple graphs, Theorem \ref{introtheorem3} $(ii)$ also gives an interesting corollary, which was first proved by Zheng \cite{Zheng2004} in the special case that $G$ was a forest.

\begin{corollary}[Corollary \ref{cor.regchordal}] \label{introcor4}
Let $G$ be a chordal graph. Then
$$\reg(\I(G)) = c + 1$$
where $c$ is the maximal number of 3-disjoint edges in $G$.
\end{corollary}

Finally, in Section 7 we show that the first syzygy module of $\I(\HH)$
when $\HH$ is properly-connected
is generated by linear syzygies
if and only if the diameter of the hypergraph $\HH$ is small enough (Theorem \ref{linearsyzygies}).  
By {\bf diameter} we mean the maximum distance between any two edges of $\HH$.
This result can be seen as the first step towards generalizing Fr\"oberg's result \cite{Fr}
characterizing graphs whose edge ideals have a linear resolution.
As an interesting corollary, if $\HH$ is a triangulated hypergraph, and 
if $\I(\HH)$ only has linear first syzygies, then
the resolution of $\I(\HH)$ must in fact be linear (Corollary \ref{linsyzcor}).


\section{Preliminaries}

We recall the relevant results concerning hypergraphs,
resolutions, and splittable ideals.

\subsection{Hypergraphs and edge ideals}  Our reference for
the hypergraph material is Berge \cite{Berge1989}.

Throughout this paper we shall assume that our hypergraphs $\HH = (\X,\E)$ are simple,
i.e., $|E| \ge 2$ for all $E \in \E$, and there is no element of $\E$ which contains another. When there is no danger of confusion, we sometimes specify a hypergraph by describing only its set of edges.

If each $E \in \E$ has the same cardinality $d$, then we call $\HH$
a {\bf $d$-uniform} hypergraph.  Note that a simple graph is a simple $2$-uniform
hypergraph.  If $\HH$ is $d$-uniform, then the associated simplicial complex
$\Delta(\HH)$ is a {\bf pure} simplicial complex, that is, all its facets have the
same dimension.

If $E$ is an edge of a hypergraph $\HH$, then we let $\HH \backslash E$ denote the
hypergraph formed by removing the edge $E$ from $\HH$.  Similarly, if $x$ is a vertex of $\HH$,
we shall write $\HH \backslash \{x\}$ to denote the hypergraph formed by removing $x$ and all
edges $E \in \E$ with the property that $x \in E$.  Note that $x$ is an isolated vertex of $\HH \backslash \{x\}$, or we can also consider
the vertex set of $\HH \backslash \{x\}$ to be $\X \backslash \{x\}$.  If $\mathcal{Y} \subset \X$,
then the {\bf induced hypergraph on $\mathcal{Y}$}, denoted $\HH_{\mathcal{Y}}$, is the
sub-hypergraph of $\HH$ whose edge set is $\{E \in \E ~|~ E \subseteq \mathcal{Y}\}$.
If there is no edge $E \in \E$ such that $E\subseteq \mathcal{Y}$, then we view
$\HH_{\mathcal{Y}}$ as the graph of the isolated vertices $\mathcal{Y}$.  

The notion of distance between edges in a hypergraph will play a fundamental
role in later discussions.  We introduce the relevant definitions here.

\begin{definition}
A {\bf chain of length $n$} in $\HH$ is a sequence
$(E_0,x_1,E_1,\ldots,x_n,E_n)$ such that
\begin{enumerate}
\item[$(1)$]  $x_1,\ldots,x_n$ are all distinct vertices of $\HH$,
\item[$(2)$]  $E_0,\ldots,E_n$ are all distinct edges of $\HH$, and
\item[$(3)$]  $x_1 \in E_0$, $x_n \in E_n$, and $x_k,x_{k+1} \in E_k$ for each $k =1,\ldots,n-1$.
\end{enumerate}
We sometimes denote the chain by $(E_0,\ldots,E_n)$ if the vertices in the chain are not being investigated. Note that $(3)$ implies that $E_i \cap E_{i+1} \neq \emptyset$ for
$i = 0,\ldots,n-1$.
If $E$ and $E'$ are two edges, then $E$ and $E'$ are {\bf connected} if there
exists a chain $(E_0,\ldots,E_n)$ where $E = E_0$ and $E' = E_n$.
If $|E| \geq |E'|$,
then the chain connecting $E$ to $E'$ is a {\bf proper chain} if
$|E_i \cap E_{i+1}| = |E_{i+1}|-1$ for all $i = 0,\ldots,n-1$.
The (proper) chain is an {\bf (proper) irredundant chain} of length $n$
if no proper subsequence is a (proper) chain from $E$ to $E'$.
\end{definition}

\begin{definition}
If $E$ and $E'$ are two edges of a hypergraph $\HH$  with $|E| \geq |E'|$, then we define the
{\bf distance} between $E$
and $E'$, denoted by $\dist_{\HH}(E,E')$, to be
\[\dist_{\HH}(E,E') = \min\{\ell ~|~ (E= E_0,\ldots,E_{\ell}=E') ~\mbox{is a proper
irredundant chain}\}.\]
If no proper irredundant chain between the two edges exists, we
set $\dist_{\HH}(E,E') = \infty$.
\end{definition}

As in the introduction, the {\bf edge ideal} of  $\HH = (\X,\E)$
is the squarefree monomial ideal
\[\I(\HH) = \left(\left\{\left. x^E = \prod_{x \in E} x ~\right|~ E \in \E\right\} \right)
\subseteq R = k[x_1,\ldots,x_n].\]
We often abuse notation and write $x^E$ for both
the edge $E$ and the corresponding monomial.

\subsection{Resolutions and splittable ideals}

Let $M$ be a graded $R$-module where  $R = k[x_1,\ldots,x_n]$.
Associated to $M$ is a {\bf minimal
graded free resolution} of the form
\[
0 \rightarrow \bigoplus_j R(-j)^{\beta_{l,j}(M)}
\rightarrow \bigoplus_j R(-j)^{\beta_{l-1,j}(M)}
\rightarrow \cdots
\rightarrow  \bigoplus_j R(-j)^{\beta_{0,j}(M)}
\rightarrow M \rightarrow 0
\]
where $l \leq n$ and $R(-j)$ is the $R$-module obtained by shifting
the degrees of $R$ by $j$.  The number $\beta_{i,j}(M)$, the $ij$th
{\bf graded Betti number} of $M$, equals the number of minimal generators
of degree $j$ in the $i$th syzygy module of $M$.

Of particular interest are the following invariants which
measure the ``size'' of the  minimal graded free resolution of $I$.
The {\bf regularity} of $I$,  denoted $\reg(I)$, is defined
by
\[\reg(I) := \max\{j-i ~|~ \beta_{i,j}(I) \neq 0\}.\]
The {\bf projective dimension} of $I$, denoted $\pdim(I)$, is defined to be
\[\pdim(I):= \max\{i ~|~ \beta_{i,j}(I) \neq 0 \}.\]
An ideal $I$ generated
by elements of degree $d$  is said to have a {\bf linear resolution}
if $\beta_{i,j}(I) = 0$ for all $j \neq i+d$.

We now recall some results concerning splittable ideals.  We use
$\G(I)$ to denote the unique minimal set of generators of a monomial
ideal $I$.

\begin{definition}[see \cite{EliahouKervaire1990}]\label{defn: split}
A monomial ideal $I$ is {\bf splittable} if $I$ is the sum
of two nonzero monomial ideals $J$ and $K$, that is, $I = J+K$, such
that
\begin{enumerate}
\item $\G(I)$ is the disjoint union of $\G(J)$ and $\G(K)$.
\item there is a {\bf splitting function}
\begin{eqnarray*}
\G(J\cap K) &\rightarrow &\G(J) \times \G(K) \\
w & \mapsto & (\phi(w),\psi(w))
\end{eqnarray*}
satisfying
\begin{enumerate}
\item for all $w \in \G(J \cap K), ~~ w = \lcm(\phi(w),\psi(w))$.
\item for every subset $S \subset \G(J \cap K)$, both
$\lcm(\phi(S))$ and $\lcm(\psi(S))$
strictly divide $\lcm(S)$.
\end{enumerate}
\end{enumerate}
If $J$ and $K$ satisfy the above properties, then
we shall say $I = J + K$ is a {\bf splitting} of $I$.
\end{definition}

When $I = J + K$ is a splitting, then there is a relation between  $\beta_{i,j}(I)$
and the graded Betti numbers of the ``smaller'' ideals.
This relation was first observed for the total Betti numbers by Eliahou and
Kervaire \cite{EliahouKervaire1990} and extended to the graded case by
Fatabbi \cite{Fatabbi2001}.

\begin{theorem}
\label{prop: ekf}
Suppose $I$ is a
splittable monomial ideal with splitting $I = J+K$.  Then
\[\beta_{i,j}(I) = \beta_{i,j}(J) + \beta_{i,j}(K) +
\beta_{i-1,j}(J\cap K) ~~\mbox{for all $i, j \geq 0$} \]
where $\beta_{i-1,j}(J \cap K) = 0$ if $i = 0$.
\end{theorem}

When $I$ is a splittable ideal, Theorem \ref{prop: ekf} gives us
the following corollary.

\begin{corollary} \label{reg pdim}
If $I$ is a splittable monomial ideal with splitting $I = J+K$, then
\begin{enumerate}
\item[$(i)$] $\reg(I) = \max\{\reg(J),\reg(K),\reg(J\cap K) - 1\}$.
\item[$(ii)$] $\pdim(I) = \max\{\pdim(J),\pdim(K),\pdim(J\cap K) + 1\}$.
\end{enumerate}
\end{corollary}

Our goal is to study the numbers $\beta_{i,j}(\I(\HH))$.
It follows directly from the definition of $\I(\HH)$
that $\beta_{0,j}(\I(\HH))$ is simply the number of edges $E \in \HH$ with $|E| = j$.
We can therefore restrict to investigating the
numbers $\beta_{i,j}(\I(\HH))$ with $i \geq 1$.  When
$\HH$ is a  $d$-uniform hypergraph, the following result
implies that we only need to consider a finite range of values of $j$ for each $i$.

\begin{theorem}
Suppose that $\HH$ is a $d$-uniform hypergraph.  If $\beta_{i,j}(\I(\HH)) \neq 0$,
then $i+d \leq j \leq \min\{n,d(i+1)\}$.
\end{theorem}

\begin{proof}
Because $\HH$ is a $d$-uniform hypergraph, $\I(\HH)$ is generated by monomials
of degree $d$.  So, $\beta_{i,j}(\I(\HH)) = 0$ for $j < i+d$, thus giving us
the lower bound.  For the upper bound, the Taylor resolution implies that
$\beta_{i,j}(\I(\HH)) = 0$ if $j > d(i+1)$.  On the other hand, Hochster's
formula implies that $\beta_{i,j}(\I(\HH)) = 0$ if $j > n$.  The conclusion
now follows.
\end{proof}


\section{Splitting edges}

Let $I$ be any squarefree monomial ideal, and suppose that $\HH$ is the
hypergraph associated to $I$, i.e., $I = \I(\HH)$.  We would like to
find splittings of $I$ so that we can make use of Theorem \ref{prop: ekf}.
In this section we describe one possible splitting of $\I(\HH)$.

One of the simplest ways to partition $\G(I)$ is to pick any $m \in \G(I)$,
and set $\G(J) = \{m\}$ and $\G(K) = \G(I) \backslash \{m\}$.  Note
that this is equivalent to picking any edge
$E$ of $\HH$, and setting
\[J = (x^E) ~~\text{and} ~~ K = \I(\HH\backslash E).\]
It is immediate that $I = \I(\HH) = J + K$, and furthermore, $J$ and $K$
satisfy condition $(1)$ of Definition \ref{defn: split}.  However, for an arbitrary
edge $E$, $J$ and $K$ may fail to satisfy condition $(2)$ of Definition \ref{defn: split}.
If $E$ is chosen so that $J$ and $K$ satisfy this condition, then
we give this edge the following name.

\begin{definition}  Let $\HH$ be a hypergraph.
An edge $E$ is a {\bf splitting edge} of $\HH$ if
\[\I(\HH) = (x^E) + \I(\HH\backslash E)\] is a splitting of $\I(\HH)$.
\end{definition}

To make use of Theorem \ref{prop: ekf}, one would therefore like a means to identify
the splitting edges of a hypergraph.  The main result of this section is the following theorem which
provides a classification of the splitting edges of a hypergraph.
This theorem answers Question 5.4.2 of \cite{HaVanTuyl2006}
which asked the equivalent question of what facet could be a splitting facet of simplicial complex.

\begin{theorem} \label{characterize: splitting facets}
Let $\HH$ be a hypergraph
with two or more edges.
Then an edge $E$ is a splitting edge of $\HH$ if and only if there exists
a vertex $z \in E$ such that
\[(x^E) \cap \I(\HH\backslash E) \subseteq (x^E) \cap \I(\HH\backslash \{z\}).\]
\end{theorem}

\begin{proof} Let $E$ be an edge of $\HH$, and
set $J = (x^E)$ and $K = \I(\HH\backslash E)$.
To prove the ``only if'' direction, we prove the contrapositive.
So, suppose that for every vertex $z \in E$, we have
\[(x^E) \cap \I(\HH\backslash E) \not\subseteq (x^E) \cap \I(\HH\backslash \{z\}).\]
Thus, for each $z \in E$, there exists a minimal generator $x^{L_z}$ of $J \cap K$
such that $x^{L_z} \not\in (x^E) \cap \I(\HH\backslash \{z\})$.
Set $S = \{x^{L_z} ~|~ z \in E\} \subseteq \G(J \cap K)$.

We will now show that no splitting function can exist.  Suppose there was a splitting
function $s:\G(J \cap K) \rightarrow \G(J) \times
\G(K)$ given by $s(w) = (\phi(w),\varphi(w)).$
Then, since $J = (x^E)$, for each $x^{L_z} \in S$, we have $\phi(x^{L_z}) = x^E$.
For each $z \in E$, let $x^{G_z} = \varphi(x^{L_z}) \in \G(K)$.  So
$G_z$ is an edge of $\HH$, and $\lcm(x^{E},x^{G_z}) = x^{E \cup G_z} =x^{L_z}$.

We claim that for each $z \in E$, we have $z \in G_z$.  Indeed,
if $z' \not\in G_{z'}$ for some $z' \in E$, then
$G_{z'}$ is an edge of $\HH \backslash \{z'\}$.  But
then $x^{L_{z'}} = \lcm(x^E,x^{G_{z'}}) = x^{E\cup G_{z'}}$
is an element of $(x^E) \cap \I(\HH\backslash \{z'\})$,
a contradiction to the choice of $x^{L_{z'}}$.

Now, since $z \in G_z$ for each $z \in E$, we have
\begin{eqnarray*}
\lcm(\varphi(S)) &=& \lcm(\{x^{G_z} ~|~ z \in E\}) = x^{\cup_{z \in E} G_z} = x^{(\cup_{z \in E} G_z) \cup E} = x^{\cup_{z \in E} (G_z \cup E)} \\
&=& x^{\cup_{z \in E} L_z} = \lcm(\{x^{L_z} ~|~ z \in E\}) = \lcm(S).
\end{eqnarray*}
But this contradicts the fact that we have a
splitting function.  This proves the ``only if'' direction.

Conversely, suppose that there exists a vertex $z$ of $E$ such that
\[(x^E) \cap \I(\HH\backslash E) \subseteq (x^E) \cap \I(\HH\backslash \{z\}).\]
This implies that $\G(J \cap K) \subseteq
\{x^{E \cup H} ~|~ H \in \HH\backslash \{z\}\}$.
We will construct a splitting function
$s = (\phi,\varphi): \G(J \cap K) \rightarrow \G(J) \times \G(K)$ which satisfies
the
conditions of Definition \ref{defn: split}.
For any $x^L \in \G(J \cap K)$ we define $\phi(x^L) = x^E \in \G(J)$. For each $x^L \in
\G(J \cap K)$, $\varphi(x^L)$ is defined as follows:
by our hypothesis, we have $L \in \{E \cup H ~|~ H \in \HH\backslash \{z\}\}$.
Thus, $\A = \{H \in \HH\backslash \{z\} ~|~ L = E \cup H\}$ is not the empty set.
We consider $\X$ as a set of alphabets (in some order of its elements) and identify each
element of $\A$ with the word formed by its
vertices (in increasing order). Let $G_L$ be the unique maximal element of $\A$
with respect to the lexicographic word ordering (which is a total order). Observe that,
by construction, $z \not\in G_L$ and $E \cup G_L = L$. Define $\varphi(x^L) = x^{G_L}$.

It is easy to see that $s = (\phi,\varphi)$ is a well defined function on $\G(J\cap K)$
and that condition (a) of Definition \ref{defn: split} is satisfied. To show that
condition (b) of Definition \ref{defn: split} is satisfied, we observe that for any
$x^L \in \G(J \cap K)$, by construction, $z$ does not divide $\varphi(x^L)$. Observe further
that for any subset $S \subseteq \G(J \cap K)$, $z$ divides $x^E$ which strictly divides
$\lcm(S)$. Thus, since $\lcm(\phi(S)) = x^E$ and since $z$ does not divide $\lcm(\varphi(S))$, we must
have that $\lcm(\phi(S))$ and $\lcm(\varphi(S))$ both strictly divide $\lcm(S)$.
The ``if'' direction is proved.
\end{proof}

\begin{remark}
Theorem \ref{characterize: splitting facets} could be reinterpreted
as describing when a squarefree monomial ideal $I = (m_1,\ldots,m_s)$
in $R = k[x_1,\ldots,x_n]$
has a splitting $I = (m_i) + (m_1,\ldots,\hat{m}_i,\ldots,m_s)$ for some
$i$.  Precisely,  $I = (m_i) + (m_1,\ldots,\hat{m}_i,\ldots,m_s)$
is a splitting
if and only if there exists a variable $x_j$ such that $x_j|m_i$ and
$(m_i) \cap (m_1,\ldots,\hat{m}_i,\ldots,m_s) \subseteq (m_i) \cap I'R$
where by $I'R$ we mean the ideal $I' = I \cap k[x_1,\ldots,\hat{x}_j,\ldots,x_n]$,
but viewed as ideal of $R$.   The result
follows from the fact that $\I(\HH\backslash \{x_j\}) =I'R$.  This reformulation
nicely illustrates that in some cases the hypergraph point of view is conceptually
easier (at least to us) to grasp.  
\end{remark}

\begin{example}
The following example illustrates that a hypergraph may not have a splitting
edge.  Let $\HH$ be the hypergraph on vertex set $\X = \{a,b,c,d,e\}$
with  edge set $\E = \{abe,ade,bce,cde\}$.  The edge ideal is
then $\I(\HH) = (abe, ade, bce, cde).$
By symmetry it suffices to show that any one of the edges is not a splitting edge.
So, consider the edge $E = abe$.  Then
\[(x^E) \cap \I(\HH \backslash E) = (abde, abce,abcde) = (abde, abce)\]
while
\[(x^E) \cap \I(\HH \backslash\{a\}) = (abce), \hspace{.25cm}
(x^E) \cap \I(\HH \backslash\{b\}) = ( abde ), ~~\text{and}~~
(x^E) \cap \I(\HH \backslash\{e\}) = (0).\]
Thus, there is no vertex $z \in E$ with the property that
$(x^E) \cap \I(\HH\backslash E) \subseteq (x^E) \cap \I(\HH\backslash \{z\})$.
\end{example}

There is a nice class of edges of a simple hypergraph that are easy to identify and
also have the property that they are splitting edges. We now define this class.

\begin{definition} \label{defn: v-leaf}
Let $\HH$ be a simple hypergraph.  An edge $E$ is a {\bf $v$-leaf} if $E$ contains
a free vertex, that is, $E$ contains a vertex $v \in \X$ such that $v$ does not belong to any
other edge of $\HH$.
\end{definition}

\begin{remark} If $\HH = G$ is a simple graph, then $v$-leaves are precisely the leaves in the usual sense.
\end{remark}

\begin{corollary} \label{splitting v-leaf}
Suppose $E$ is a $v$-leaf of a hypergraph $\HH$.  Then $E$
is a splitting edge of $\HH$.
\end{corollary}

\begin{proof} If $v$ is the free vertex in $E$, then $\HH\backslash E = \HH\backslash\{v\}$.
Now apply Theorem \ref{characterize: splitting facets}. 
\end{proof}

S. Faridi \cite{Faridi2002} introduced the notion of a leaf for a simplicial
complex $\Delta$.  Precisely, a facet $F$ of $\Delta$ is a {\bf leaf} if $F$ is the only
facet of $\Delta$, or there exists a facet $G \neq F$ in $\Delta$ such that
$F \cap F' \subseteq F \cap G$ for all facets $F' \neq F$ in $\Delta$.
We can translate Faridi's definition into hypergraph language; we call the
translated version of Faridi's leaf a $f$-leaf to distinguish it from a $v$-leaf.

\begin{definition} An edge $E$ of a hypergraph $\HH$ is a {\bf $f$-leaf} if $E$
is the only edge of $\HH$, or if there exists an edge $H$ of $\HH$ such that
$E \cap E' \subseteq E \cap H$ for all edges $E' \neq E$ of $\HH$.
\end{definition}

We introduce two types of hypertrees and hyperforests based upon the two notions of leaves.

\begin{definition}
A hypergraph $\HH$ is a {\bf $v$-forest}, respectively, $f$-{\bf forest}, if every
induced subgraph of $\HH$, including $\HH$ itself, contains a $v$-leaf, respectively,
a $f$-leaf.  If $\HH$ is connected, we call $\HH$ a {\bf $v$-tree}, respectively,
{\bf $f$-tree}.  When $\HH$ is a $f$-forest, the associated simplicial
complex $\Delta(\HH)$ is called a {\bf simplicial forest}.
\end{definition}

Notice that when $\HH = G$ is a simple graph, the notions of $v$-leaf and $f$-leaf coincide.
So, with simple graphs, the notions of a $v$-forest and a $f$-forest coincide with the usual 
notion of a forest.
These definitions, however,  are not equivalent in a general hypergraph, as illustrated below.

\begin{example}\label{ex: v-leaf}
A $f$-leaf must always contain a free vertex (cf. \cite[Remark 2.3]{Faridi2002}), thus
every $f$-leaf is a $v$-leaf.  However, a $v$-leaf need not be a $f$-leaf.  For example,
consider the hypergraph $\HH$ on $\X = \{a,b,c,d,e,f\}$ with the edge set
$\E = \{abf, bcd, def\} = \{E_1,E_2,E_3\}$.
Each edge is a $v$-leaf since each edge has a vertex not in the other two edges.
However, $\HH$ has no $f$-leaf. By symmetry, it is enough to show that $E_1 = abf$ cannot
be a $f$-leaf.  Indeed, $E_1 \cap E_2 \not\subseteq E_1 \cap E_3$ and $E_1 \cap E_3 \not\subseteq
E_1 \cap E_2$.

The hypergraph $\HH$ is an example of a $v$-tree, but $\HH$ is not a $f$-tree since
$\HH$ has no $f$-leaf, although all its induced subgraphs have a $f$-leaf.
\end{example}

Because a $f$-leaf is a $v$-leaf, Corollary \ref{splitting v-leaf}
immediately gives:

\begin{corollary} \label{splitting leaf}
If $E$ is a $f$-leaf of a hypergraph $\HH$, then $E$ is a splitting edge of $\HH$.
\end{corollary}


\section{Properly-connected hypergraphs}

Given a hypergraph $\HH$, we would like to express the numbers $\beta_{i,j}(\I(\HH))$
in terms of the graded Betti numbers of edge ideals associated to subgraphs of $\HH$;
this would lead to recursive-type formulas.
When $E$ is a splitting edge of a hypergraph $\HH$, Theorem
\ref{prop: ekf} implies that $\beta_{i,j}(\I(\HH))$ can
be computed from the graded Betti numbers of the ideals $(x^E)$, $\I(\HH\backslash E)$,
and $L = (x^E) \cap \I(\HH\backslash E)$.  The Betti numbers of $(x^E)$ are trivial
to compute, while those of $\I(\HH\backslash E)$ already
correspond to the edge ideal of a sub-hypergraph of $\HH$. Thus one only
needs to relate the numbers $\beta_{i,j}(L)$
to the Betti numbers of an edge ideal of some other sub-hypergraph.  For a general hypergraph,
this appears to be a difficult problem.

The goal of this section is to introduce a family of $d$-uniform hypergraphs,
which we call properly-connected, that among other things enables us  to relate the
graded Betti numbers of $L$ to those of an edge ideal associated to a sub-hypergraph 
of $\HH$.

\begin{definition}
A $d$-uniform hypergraph $\HH = (\X,\E)$ is said to be {\bf properly-connected} if
for any two edges $E$ and $E'$ of
$\HH$ with the property that $E \cap E' \neq \emptyset$,
then
\[\dist_{\HH}(E,E') = d - |E\cap E'|.\]
Otherwise, we say $\HH$
is {\bf not properly-connected}. 
\end{definition}

\begin{remark}
Our definition of properly-connected is similar to (but not equivalent to) what
Zheng \cite[Definition 3.14]{Zheng2004} called the {\bf intersection property}
for a simplicial complex.  If $\Delta$ is a pure simplicial forest, then $\Delta$ has
the intersection property if for any two facets $F,F' \in \Delta$
the distance between $F$ and $F'$ (defined in terms
of the length of chain of between the two facets) is determined by $|F \cap F'|$.
\end{remark}

\begin{example}
Consider the $4$-uniform hypergraph $\HH$ with edge set
$$\E = \{x_1x_2x_3x_4, x_1x_2x_3x_7, x_1x_2x_6x_7, x_1x_5x_6x_7,
x_1x_5x_6x_8\}.$$
There is a proper irredundant chain of length $4$ from the
edge $E = x_1x_2x_3x_4$ to $E' = x_1x_5x_6x_8$ (to
form the chain, just take the edges as listed in $\E$).  Furthermore,
there is no shorter such chain.  But $E$ and $E'$ have
a nonempty intersection.  So $\HH$ is not properly-connected
since  $4 = \dist_{\HH}(E,E') \neq 4 - |E \cap E'| = 3$.
This hypergraph is not properly-connected.
\end{example}

\begin{example} Every finite simple graph $G$ is properly-connected.
To see this, note that a graph is clearly a $2$-uniform hypergraph.
If $E,E'$ are two edges of $G$ such that $E \cap E' \neq \emptyset$, then
either $E$ and $E'$ are the same edge,
or $E$ and $E'$ share exactly one vertex.  In the first case,
$\dist_{G}(E,E') = 2 - |E \cap E'| = 2-2 =0$, while in the second case
$\dist_{G}(E,E') = 2 - |E \cap E'| = 1$.  
So, in this sense, properly-connected hypergraphs generalize simple graphs.
\end{example}

Properly-connected hypergraphs are appealing combinatorial objects to study
because within this family, the notions of $v$-leaf and $f$-leaf become equivalent.
As well, splitting edges of properly-connected hypergraphs can be described
combinatorially.  We prove both of these assertions.

\begin{theorem}
Suppose $\HH$ is a $d$-uniform properly-connected hypergraph, and $E$ is an edge of $\HH$.  Then
$E$ is a $v$-leaf if and only if $E$ is a $f$-leaf.
\end{theorem}

\begin{proof} 
Because we know an $f$-leaf is a $v$-leaf, it suffices to prove the converse.
If $\HH$ has only one edge then we are done.  So, suppose that $\HH$ has
at least two edges.  Let $E$ be a $v$-leaf with free vertex $v$.
Let $H$ be any edge of $\HH$ with $H \cap E \neq \emptyset$. 
If there is no such $H$, then
$E$ is automatically a $f$-leaf.   
Since $\HH$ is properly-connected, there is a proper
chain $E_0 = E, E_1,\ldots, E_k =H$ from $E$ to $H.$  Because $|E| = |E_1| = d$ and $|E \cap E_1| =
d-1,$ $E \cap E_1 = E \backslash \{v\}.$
To see that $E$ is a $f$-leaf, let $G$ be any other edge of $\HH.$
Then $E \cap G \subset E\backslash\{v\} = E \cap E_1.$
\end{proof}

Let $E$ be an edge of a $d$-uniform properly-connected hypergraph $\HH$.  If $H$ is
any edge of $\HH$ with $\dist_{\HH}(E,H)=1$, then $|H \backslash E| = 1$, or
in other words, $H \backslash E = \{z\}$ for some vertex $z$.
Before classifying splitting edges, we introduce the following definition.

\begin{definition}
If $E$ is an edge of a $d$-uniform properly-connected hypergraph $\HH$, then 
the {\bf vertex neighbor set of $E$} is the following subset of $\X$:
\[N(E) =  \bigcup_{\{H \in \HH ~|~ \dist_{\HH}(E,H) =1\}} H \backslash E.\]
\end{definition}

\begin{example}
When $G$ is a finite simple graph, and $x$ is a vertex, then $N(x)$ denotes
all the neighbors of $x$.  If $E = \{u,v\}$ is any edge of $G$, then
$N(E) = (N(u) \cup N(v)) \backslash \{u,v\}$.
\end{example}

\begin{theorem}\label{char: split edge}
Let $E$ be an edge of a $d$-uniform properly-connected hypergraph $\HH$,
and suppose $N(E) = \{z_1,\ldots,z_t\}$.  Then
$E$ is a splitting edge if and only if there exists a vertex $z \in E$
such that $(E \backslash \{z\}) \cup \{z_i\} \in \HH$ for each
$z_i \in N(E)$.
\end{theorem}

The proof of this theorem depends upon the following two lemmas.

\begin{lemma}\label{lem: chains}
Let $\HH$ be a $d$-uniform properly-connected hypergraph.
Suppose that $E= E_0 = \{x_1, \ldots, x_d \}$ and $E'$ are edges in $\HH$ with
$\dist_{\HH}(E,E') = t \leq d$.  Then, after relabelling,
there exist edges $E_1, \ldots, E_t$ such that
$E_i = \{y_1, \ldots, y_i, x_{i+1}, \ldots, x_d\},
$ $E_t = E',$ and $y_i \notin E_{j}$ for all $j < i.$
\end{lemma}

\begin{proof}
Since $\dist_{\HH}(E,E') = t$, there must be a proper irredundant chain of
edges $E_0 = E, \ldots, E_t = E'$.  Since $E_i$ differs from $E_{i+1}$ by exactly one vertex,
for each $i,$ $|E \cap E_i| \geq d-i$ because at most one vertex changes at each stage.
Since $(E_0, \dots, E_t)$ is an irredundant chain and $\HH$ is properly-connected, for $i < d$, we must have
\[i = \dist_{\HH}(E_0, E_i) = d - |E_0 \cap E_i|.\]
Hence,
$|E_0 \cap E_i| = d-i$ for any $i$ less than $d$ for which the expression makes sense.
Moreover, if $i = t = d$, then $\dist_\HH(E_0,E_i) = d,$ and we have $E_0 \cap E_i =
\emptyset$. That is, $|E_0 \cap E_i| = 0 = d-i$.

We will prove the result using induction on $i$.
Let $E =E_0 = \{x_1,\ldots,x_d\}$, and assume the vertices are labeled so that
$x_1 \notin E_1$.  We know that $|E_0 \cap E_1| = d-1$ which implies that $E_1 =
\{y_1, x_2, \ldots, x_d\}$ where $y_1 \notin E_0$, thus
proving the base case.

Now assume that $E_0, \ldots, E_i$ satisfy the claim, i.e., that
$E_i = \{y_1, \ldots, y_i, x_{i+1}, \ldots, x_d\}$ with $y_i \notin E_j$ for all $j < i$.
We know that $|E_i \cap E_{i+1}| = d-1,$ so that $E_{i+1}$ is constructed from $E_i$ by
removing some vertex and adding a vertex that we will call $y_{i+1}$ which is not in $E_i.$
First, we claim that the vertex that we remove from $E_i$ cannot be a $y_j.$  If we were to
replace some $y_j$ with a vertex $y_{i+1},$ then $|E_0 \cap E_i| = d-i \leq |E_0 \cap E_{i+1}|$
which contradicts our earlier assumption that $|E_0 \cap E_{i+1}| = d-i-1$.
So, we may assume that $y_{i+1}$ replaces $x_{i+1}.$  If, $y_{i+1} = x_j$ for some $j \le i$ then
$|E_0 \cap E_{i+1}| =
|E_0 \cap E_i|$ which is a contradiction as before.  Therefore,
$y_{i+1} \notin E_j$ for any $j \leq i.$
\end{proof}

\begin{lemma}\label{intersectionideal}
Let $E$ be any edge of a $d$-uniform 
properly-connected hypergraph $\HH$.  Then
\begin{eqnarray*}
(x^E) \cap \I(\HH\backslash E)& =& (\{\lcm(x^E,x^H) ~|~ H \in \HH ~\mbox{and}~
\dist_{\HH}(E,H) =1\}) + \\
&&(\{\lcm(x^E,x^H) ~|~ H \in \HH ~\mbox{and}~\dist_{\HH}(E,H) \geq d+1\}).
\end{eqnarray*}
\end{lemma}

\begin{proof}  Set
\begin{eqnarray*}
A & =& (\{\lcm(x^E,x^H) ~|~ H \in \HH \backslash E ~\mbox{and}~\dist_{\HH}(E,H) \leq d \}) ~~\mbox{and}\\
B & =&(\{\lcm(x^E,x^H) ~|~ H \in \HH ~\mbox{and}~\dist_{\HH}(E,H) \geq d+1\}).
\end{eqnarray*}
By definition $(x^E) \cap \I(\HH\backslash E) = A + B$.  Thus, if
we set
\[C = (\{\lcm(x^E,x^H) ~|~ H \in \HH ~\mbox{and}~\dist_{\HH}(E,H) =1\}),\] then 
it suffices to show that
$A =  C$.  Since $C \subseteq A$ is clear, we now show the reverse containment.

Let $x^{E \cup H} = \lcm(x^E,x^H)$ be a generator of $A$, i.e., suppose $H \in \HH \backslash E$ and
$t = \dist_{\HH}(E,H) \leq d$.  Note that we can assume that $2 \leq t \leq d$
because if $t= \dist_{\HH}(E,H) = 1$, then $x^{E \cup H} \in C$.  So there exists
a proper irredundant chain $E = H_0,H_1,H_2,\ldots,H_t = H$ whose length is minimal
among all proper irredundant chains from $E$ to $H$.


Now if $E = \{x_1,\ldots,x_d\}$,
then $H_1 = \{x_1,\ldots,\hat{x}_i,\ldots,x_d,z\}$ where by $\hat{x}_i$ we mean
the
vertex $x_i$ is removed, and $z$ is not one of $x_1,\ldots,x_d$.
From this observation, we have
\[\lcm(x^E,x^{H_1}) = x^{E \cup \{z\}} = x^Ez.\]
Now $x^Ez$ is a generator of $C$.
To finish the proof,  Lemma \ref{lem: chains} implies that $z \in H_i$
for $i = 2,\ldots,t$.  Therefore, $\lcm(x^E,x^{H_i}) = x^{E \cup H_i}$ is divisible
by $x^Ez$, and thus is in $C$.  In particular $x^{E \cup H} \in C$.
\end{proof}

\noindent{\it Proof of Theorem \ref{char: split edge}.}
Suppose that $E$ is a splitting edge.  By Theorem \ref{characterize: splitting facets}
there is a vertex $z \in E$ such that
$(x^E) \cap \I(H\backslash E) \subseteq
(x^E) \cap \I(H\backslash\{z\})$.
Let $z_i \in N(E)$.  We will show that $(E \backslash \{z\}) \cup \{z_i\}$ is an edge of
$\HH\backslash \{z\} \subseteq \HH$.  Since $z_i \in N(E)$, there exists an edge $H$
with $\dist_{\HH}(E,H) = 1$ such that $H \backslash E = \{z_i\}$.
Thus, $x^{E \cup H}$ is a generator of $(x^E) \cap \I(\HH\backslash E)$.
We thus must have $x^{E \cup H} \in (x^E) \cap \I(\HH\backslash\{z\})$.
Hence there is an edge $H' \in \HH\backslash \{z\}$ such that $E \cup H = E \cup H'$. Because
$|E \cap H| = d-1$,  we must have that $|E \cap H'| = d-1$.
Since $z \not\in H'$ and $z_i \not\in E$, we must have
$H' = (E \backslash \{z\}) \cup \{z_i\}$.  So, $(E \backslash \{z\}) \cup \{z_i\}
\in \HH\backslash\{z\}$ as desired.

Conversely, suppose there exists a vertex $z \in E$ such
that $(E \backslash \{z\}) \cup \{z_i\} \in \HH$ for each
$z_i \in N(E)$.  Let $x^L$ be any minimal generator of $(x^E) \cap \I(\HH\backslash E)$.
By Lemma \ref{intersectionideal}, we have $L = E \cup H$ with $\dist_{\HH}(E,H) = 1$
or $L = E \cup H$ with $\dist_{\HH}(E,H) \geq d+1$.  If $\dist_{\HH}(E,H) \geq d+1$,
then $z \not\in H$ since $E \cap H = \emptyset$.  So $H \in \HH\backslash \{z\}$,
and hence $x^L \in (x^E) \cap \I(\HH\backslash \{z\})$.  So,
suppose $L = E \cup H$ with $\dist_{\HH}(E,H)=1$.  Then, the exists
$z_i \in N(E)$ such that $E \cup H = E \cup \{z_i\}$.  By our hypothesis, the
edge $E' = (E \backslash \{z\}) \cup \{z_i\} \in  \HH$.  But then $E' \in \HH\backslash \{z\}$.
Furthermore, $L = E \cup H = E \cup E'$.  So $x^L \in (x^E) \cap \I(\HH\backslash \{z\})$.
We have now shown that $(x^E) \cap \I(\HH\backslash E) \subseteq
(x^E) \cap \I(\HH\backslash \{z\})$, so by Theorem \ref{characterize: splitting facets}
the edge $E$ must be a splitting edge.
\qed

\begin{example} 
We give an example of a $3$-uniform properly-connected hypergraph
which has a splitting edge that is not a $v$-leaf.  Let $\HH$ be the hypergraph
with edge set
\[\E = \{x_1x_2x_3,x_1x_2x_4,x_1x_3x_5,x_2x_3x_4,x_2x_3x_5,x_3x_4x_5\}.\]
One can verify that $\HH$ is properly connected by showing that
$\dist_{\HH}(E,E') = 3 - |E \cap E'|$ for every pair of edges in $\E$.
Now $E = x_1x_2x_3$ is not a $v$-leaf because it does not contain a free vertex. 
We can use Theorem \ref{char: split edge} to verify that $E$ is a splitting edge.
In this case $N(E) = \{x_4,x_5\}$ since the edges of distance one from
$E$ are $\{x_1x_2x_4,x_1x_3x_5,x_2x_3x_4,x_2x_3x_5\}$.  Then $E$
is a splitting edge since $(E\backslash\{x_1\}) \cup \{x_4\} = x_2x_3x_4$ and $(E\backslash\{x_1\}) \cup \{x_5\} = x_2x_3x_5$ are both edges of $\HH$.  Note that even when
$E$ is a splitting edge, the graph $\HH\backslash E$ may fail to
be properly-connected.  In this case, if we remove $E$ from $\HH$, the resulting
hypergraph fails to be properly-connected because edges $E_1= x_1x_2x_4$ and $E_2= x_1x_3x_5$
intersect at $x_1$, but there is no proper chain of length $2 = 3 -|E_1 \cap E_2|$
in $\HH\backslash E$ between these two edges.
\end{example}

\begin{notation} \label{H'}
Suppose $E$ is an edge of a $d$-uniform properly-connected hypergraph $\HH$. 
For simplicity of notation, 
throughout the rest of the paper, when not specified, 
$\HH'$ refers to the sub-hypergraph
\[\HH' = \{H \in \HH ~|~ \dist_\HH(E,H) \ge d+1\}.\] 
\end{notation}

The following lemma tells us the properly-connected property is passed on to $\HH'$.

\begin{lemma}\label{H'prop-con}
If $E$ is an edge of a $d$-uniform properly-connected hypergraph $\HH$, then
$\HH'$ is also a $d$-uniform properly-connected hypergraph.
\end{lemma}

\begin{proof}
Because it is clear that $\HH'$ is a $d$-uniform hypergraph, it suffices
to show that $\HH'$ is properly connected.  So, suppose that the edges $H,H' \in \HH'$
have the property that $H \cap H' \neq \emptyset$.  Because they are also edges of $\HH$,
there exists a chain $H = H_0,H_1,\ldots,H_t =H'$ in $\HH$
such that $t = \dist_{\HH}(H,H') = d - |H\cap H'|$.  If all the
edges $H_i$ for $i=1,\ldots,t-1$ are also in $\HH'$,
then it is clear that $t = \dist_{\HH'}(H,H') = d-|H\cap H'|$.
So, suppose there is an edge $H_i$ in the chain with $i \in \{1,\ldots,t-1\}$
and $H_i \not\in \HH'$.  Then $s = \dist_{\HH}(E,H_i) \leq d$.  
Let $E=E_0,E_1,\ldots,E_s=H_i$ be the proper irredundant chain
in $\HH$ between $E$ and $H_i$.  Then $\dist_{\HH}(E_1,H_i) = s-1 < d$.
But this means that $|E_1 \cap H_i| \neq \emptyset$.  Let $x \in E_1 \cap H_i$.
By Lemma \ref{lem: chains} the vertex $x$ must be in either $H$ or $H'$.
Without loss of generality, assume $x \in H$.  But then $\dist_{\HH}(E_1,H) =
d - |H \cap E_1| \leq d-1$.  But since $E$ is distance one from $E_1$,
this means there is a proper chain of length $d$ from $E$ to $H$, contradicting
the fact that $H \in \HH'$.
\end{proof}
As a byproduct of Lemma \ref{intersectionideal}, we can rewrite
$(x^E) \cap \I(\HH\backslash E)$ in terms of the edge ideal of $\HH'$.

\begin{corollary}
Let $E$ be any edge of a $d$-uniform properly-connected hypergraph $\HH$,
and suppose $N(E) = \{z_1,\ldots,z_t\}$.  Then
\[(x^E) \cap \I(\HH \backslash E) = x^E((z_1,\ldots,z_t) + \I(\HH')).\]
\end{corollary}

\begin{proof}
It is straight forward to verify that
\[x^E(z_1,\ldots,z_t) = (\{\lcm(x^E,x^H) ~|~ H\in \HH ~\mbox{and}~\dist_{\HH}(E,H) =1\}).\]
If $H \in \HH \backslash E$ with $\dist_{\HH}(E,H) \geq d+1$, then
because $\HH$ is properly-connected, $|E \cap H| = \emptyset$.
So
\[x^E\I(\HH')=
(\{\lcm(x^E,x^H) ~|~ H \in \HH ~\mbox{and}~\dist_{\HH}(E,H) \geq d+1\}).\]
The result now follows from Lemma \ref{intersectionideal}.
\end{proof}

When $E$ is an edge of a properly-connected hypergraph, we can also
describe the graded Betti numbers of $(x^E) \cap \I(\HH\backslash E)$
in terms of those of $\I(\HH')$.

\begin{lemma} \label{betti intersection}
Let $E$ be any edge of a $d$-uniform properly-connected hypergraph $\HH$.
Set $t = |N(E)|$.  Then
\[\beta_{i-1,j}((x^E) \cap \I(\HH\backslash E)) =
 \sum_{l=0}^i \binom{t}{l}
\beta_{i-1-l,j-d-l}(\I(\HH'))\]
where $\beta_{-1,j}(\I(\HH')) = 1$ if $j = 0$ and $0$ otherwise.
\end{lemma}

\begin{proof}  If $N(E) = \{z_1,\ldots,z_t\}$,
then by the previous corollary,
\begin{eqnarray*}
\beta_{i-1,j}((x^E) \cap \I(\HH\backslash E)) &= &\beta_{i-1,j}(x^E((z_1,\ldots,z_t)+
\I(\HH'))) \\
&=& \beta_{i-1,j-d}((z_1,\ldots,z_t)+ \I(\HH'))\\
&=&
\beta_{i,j-d}(R/((z_1,\ldots,z_t)+ \I(\HH')).
\end{eqnarray*}
None of the generators of $\I(\HH')$ are divisible by $z_i$ for
$i = 1,\ldots,t$.  To see this, suppose that $x^H \in \I(\HH')$
is divisible by some $z_i$, i.e., $z_i$ is a vertex of the edge $H$.
Now there is a edge $H_i$ with $z_i \in H_i$ and $\dist_{\HH}(E,H_i)=1$.
Since $H \cap H_i \neq \emptyset$
and because
$\HH$ is properly-connected, $p =\dist_{\HH}(H,H_i) = d - |H \cap H_i| < d$.
So there is a proper irredundant chain $H_i = H'_0,\ldots,H'_p = H$.  But
then $E,H_i=H'_0,\ldots,H'_p=H$ forms a proper irredundant chain of length $p+1 \leq d$,
and thus $\dist_{\HH}(E,H) \leq d$, contradicting the fact that $\dist_{\HH}(E,H)\geq d+1$.

We modify our notation and write $R = k[z_1,\ldots,z_t,x_1,\ldots,x_s]$ where
 $\{x_1,\ldots,x_s\} = \X \backslash N(E)$.
Then
\[ R/((z_1,\ldots,z_t) + \I(\HH')) \cong R_1/(z_1,\ldots,z_t) \otimes_k
R_2/\I(\HH')\]
where $R_1 = k[z_1,\ldots,z_t]$ and $R_2 = k[x_1,\ldots,x_s]$, and
where we view $\I(\HH')$ as an ideal of $R$ and as the ideal of $R_2$
generated by the same elements.
By tensoring the resolutions of  $R_1/(z_1,\ldots,z_t)$ and $R_2/\I(\HH')$
together we get (see, for example, Lemma 2.1 and Corollary 2.2 of \cite{JK})
\small
\[\beta_{i,j-d}(R/L) =
\sum_{l_1=0}^i\sum_{l_2 =0}^{j-d} \beta_{l_1,l_2}(R_1/(z_1,\ldots,z_t))
\beta_{i-l_1,j-d-l_2}(R_2/(\I(\HH'))\]
\normalsize
where $L = (z_1,\ldots,z_t) + \I(\HH')$.
Since $z_1,\ldots,z_t$ is a regular sequence on $R_1$
\[
\beta_{l_1,l_2}(R_1/(z_1,\ldots,z_t)) =
\left\{
\begin{array}{ll}
0 & \mbox{if $l_2 \neq l_1$} \\
\binom{t}{l} & \mbox{if $l = l_2 = l_1$.}
\end{array}
\right.\]
As a consequence, the previous expression reduces to
\[\beta_{i,j-d}(R/L) = \sum_{l=0}^i \binom{t}{l}
\beta_{i-l,j-d-l}(R_2/\I(\HH')).\]
We are now done since
\[
\beta_{i-l,j-d-l}(R_2/\I(\HH')) = \beta_{i-l,j-d-l}(R/\I(\HH'))
= \beta_{i-l-1,j-d-l}(\I(\HH'))\]
for all $l$ (where we adopt the convention that
$\beta_{-1,j}(\I(\HH')) = 1$ if $j=0$ and $0$ if $j \neq 0$).
\end{proof}

When $\HH$ is a properly-connected hypergraph, 
we obtain the following recursive like formula for $\beta_{i,j}(\I(\HH))$.
This result generalizes a similar result for  simple
graphs found in \cite{HaVanTuyl2005}.

\begin{theorem} \label{theoremsplitting}
Let $\HH$ be a $d$-uniform properly-connected hypergraph and let
$E$ be a splitting edge of $\HH$. Suppose 
$\HH' = \{H \in \HH ~|~ \dist_{\HH}(E,H) \geq d+1\}$, and $t= |N(E)|$.
Then for all $i \geq 1$
\[\beta_{i,j}(\I(\HH)) = \beta_{i,j}(\I(\HH\backslash E))
+ \sum_{l=0}^i \binom{t}{l}
\beta_{i-1-l,j-d-l}(\I(\HH')).\]
Here, $\beta_{-1,j}(\I(\HH')) = 1$ if $j =0$ and $0$ if $j \neq 0$.
\end{theorem}

\begin{proof} Since $E$ is a splitting edge, by Theorem \ref{prop: ekf} we have
\[\beta_{i,j}(\I(\HH)) = \beta_{i,j}((x^E))+ \beta_{i,j}(\I(\HH\backslash E))
+ \beta_{i-1,j}((x^E) \cap \I(\HH\backslash E)).\]
When $i \geq 1$, $\beta_{i,j}((x^E)) = 0$.  Now substitute
the formula of Lemma \ref{betti intersection} into the last expression.
\end{proof}


\section{Triangulated properly-connected hypergraphs}

If $\HH$ is a properly-connected hypergraph with splitting edge $E$,
the sub-hypergraphs $\HH\backslash E$ and $\HH'$ in Theorem \ref{theoremsplitting}
may or may not have a splitting edge.  In fact, $\HH\backslash E$
may not even be a properly-connected hypergraph.
These facts prevents us from
using Theorem \ref{theoremsplitting} to recursively compute $\beta_{i,j}(\I(\HH))$
for any hypergraph.  One is lead to ask if there is any subfamily of properly-connected
hypergraphs for which the formula is recursive. In this section, we
introduce one such family which generalizes the notion of a chordal graph.  
In \cite{HaVanTuyl2005} it was shown that 
hyperforests (i.e., a simplicial forest in the sense of \cite{Faridi2002}) is a family 
of hypergraphs
for which the graded Betti numbers can be computed
recursively.  Since a hyperforest need not be properly-connected, the results
of this section give a partial generalization of \cite{HaVanTuyl2005}.

We begin by recalling the definition of a chordal graph.

\begin{definition}
A graph $G$ is called {\bf chordal} if every cycle of length 4 or larger has
a chord, that is, an edge joining two  nonadjacent vertices in the the
cycle.
\end{definition}

An alternative characterization for chordal graphs can be found in \cite{PTW} 
(due to Dirac \cite{D}). This characterization will prove more
suitable when generalizing to properly-connected hypergraphs.

\begin{theorem} \label{chordal-triangulated}
A graph $G$ is chordal if and only if every induced subgraph of $G$ contains
a vertex $v$ whose neighborhood $N(v)$ is a complete graph.
\end{theorem}

In the above theorem, because $v$ is adjacent to every vertex in $N(v)$,
we also have that the induced graph on $N(v) \cup \{v\}$ is also a complete graph.
To extend this definition, we first introduce
an analog of complete graphs.

\begin{definition}
The {\bf $d$-complete hypergraph of order $n$}, denoted by $\K_n^d$,
is the hypergraph consisting of all the $d$-subsets of 
the vertex set $\X$, where
$|\X| = n$.  When $d = 2$, then $\K_n^2$ is the usual complete
graph $\K_n$.  When $n < d$, we consider $\K_n^d$ as the
hypergraph with $n$ isolated vertices.
If $n=0$, then $\K_0^d$ is the empty graph which
we view as the the $d$-complete hypergraph of order $0$. 
\end{definition}

\begin{definition}
Two distinct vertices $x,y \in \X$ are {\bf neighbors} if there is an edge $E \in \HH$
such that $x,y \in E$.  For any vertex $x \in \X$, the
{\bf neighborhood of $x$}, denoted $N(x)$, is the set
\[N(x) = \{y \in \X ~|~ \mbox{$y$ is a neighbor of $x$}\}.\]
\end{definition}

Observe that if $E$ is any edge of $\HH$ and $x \in E$, then $E \subseteq N(x) \cup \{x\}$.

\begin{definition}
A $d$-uniform properly-connected hypergraph $\HH$ is
said to be {\bf triangulated} if for every nonempty subset $\Y \subseteq \X$, the
induced subhypergraph $\HH_\Y$ contains a vertex $x \in \Y \subseteq \X$
such that the induced hypergraph of $\HH_\Y$ on
$N(x) \cup \{x\}$ is a  $d$-complete hypergraph of order $|N(x)|+1$.
\end{definition}

By virtue of Theorem \ref{chordal-triangulated}, the simple graphs that are triangulated are precisely
the chordal graphs.
We shall show that properly-connected hyperforests are triangulated hypergraphs.

\begin{theorem}
Suppose that $\HH$ is a $d$-uniform properly-connected hypergraph that is a
$v$-forest (or equivalently, $f$-forest).  Then $\HH$ is a triangulated hypergraph.
\end{theorem}

\begin{proof}
For any $\Y \subseteq \X$, the induced subgraph $\HH_{\Y}$ must contain a $v$-leaf, say $E$.
Since $E$ is a $v$-leaf,
$E$ contains a free vertex, say $x$.  Suppose $E = \{x,x_2,\ldots,x_d\}$. Then
$N(x) = \{x_2,\ldots,x_d\}$.  But the induced graph of $\HH_{\Y}$ on $N(x) \cup\{x\}$ is the
simply the edge $E$ which is the
$d$-uniform complete hypergraph $\K_{d}^d$.  So $\HH$ is a triangulated hypergraph.
\end{proof}

The following lemma is the key result needed to prove that Theorem \ref{theoremsplitting}
is recursive for triangulated hypergraphs.

\begin{lemma} \label{triangulatedlemma}
Let $\HH$ be a triangulated hypergraph.  Then there exists
an edge $E \in \HH$ such that
\begin{enumerate}
\item[$(a)$] $E$ is a splitting edge, and
\item[$(b)$] the subgraphs $\HH\backslash E$ and
$\HH'$ are triangulated hypergraphs.
\end{enumerate}
\end{lemma}

\begin{proof}
Since $\HH$ is a triangulated hypergraph, there exists a vertex $x \in \X$ such
that the induced hypergraph on $N(x) \cup \{x\}$ is a $d$-complete hypergraph.
Let $E$ be any edge of $\HH$ that contains $x$.  We will show that $E$ is
an edge that satisfies $(a)$ and $(b)$.

$(a)$
Suppose that $N(E) = \{z_1,\ldots,z_t\}$.  For each $z_i \in N(E)$, there
must be an edge $E_i \in \HH$ such that $\dist_{\HH}(E,E_i) = 1$
and $E \cup E_i = E \cup \{z_i\}$.  For each $i$,
either $x \in E_i$ or $x \not\in E_i$.  If $x \not\in E_i$, then
$(E \backslash \{x\}) \cup \{z_i\} = E_i \in \HH$.
Now, suppose $x \in E_i$.
Since $z_i \in E_i$, we have $z_i \in N(x)$.  If $E = \{x,x_2,\ldots,x_d\}$, then
$\{x_2,\ldots,x_d,z_i\} \subseteq N(x)$ is a subset of size $d$ in $N(x) \cup \{x\}$.
But since the  induced hypergraph on $N(x) \cup \{x\}$ is a $d$-complete hypergraph,
that means that $\{x_2,\ldots,x_d,z_i\}$ is an edge of $\HH$.  This edge
is simply $(E \backslash \{x\}) \cup \{z_i\}$.  So, $E$ is a splitting edge
by Theorem \ref{char: split edge}.

$(b)$  We consider $\HH\backslash E$ first.
We begin by showing that $\HH\backslash E$ is properly-connected.
If $H,H' \in \HH\backslash E$ with $H \cap H' \neq \emptyset$,
then in $\HH$ we also have $H \cap H' \neq \emptyset$.  Since
$\HH$ is properly-connected, we can find a proper irredundant chain
$H=E_0,E_1,\ldots,E_t=H'$ where $t = \dist_{\HH}(H,H') = d-|H\cap H'|$.
If $E \not\in \{E_1,\ldots,E_{t-1}\}$, then this chain remains
a proper irredundant chain in $\HH\backslash E$ giving
us $t = \dist_{\HH\backslash E}(H,H')=d-|H\cap H'|$.

So suppose $E \in \{E_1,\ldots,E_{t-1}\}$.  Let $x \in E$ be
the vertex such that the induced hypergraph on $N(x) \cup
\{x\}$ is a $d$-complete hypergraph.  Let $E_{i-1}$ and $E_{i+1}$
be the edges that appear immediately before and after $E$, respectively,
in the chain $E_0,\ldots,E_t$.  There then exists a vertex $z_{i-1} \in E_{i-1}$
such that $\{z_{i-1}\} = E_{i-1} \setminus E$, and a vertex 
$z_{i+1} \in E_{i+1}$ such that $\{z_{i+1}\} = E_{i+1} \setminus E$.
By Lemma \ref{lem: chains} there are three cases to consider:
(i) $x \in E_{i-1},E,$ and $E_{i+1}$, (ii) $x \in E_{i-1}$ and $E$, but $x \not\in E_{i+1}$,
or (iii) $x \not\in E_{i-1}$ but $x \in E$ and $E_{i+1}$ (Lemma \ref{lem: chains}
shows that when moving through the chain, one removes one vertex from
an edge and replaces it with another vertex, and furthermore,
once you add a vertex to a chain, this vertex appears in all later edges
in the chain.) In case (i), let $E' = E_{i-1} \cap E_{i+1}$. Note
that $|E'|=d-2$. Then $E' \cup \{z_{i-1},z_{i+1}\}$ is a subset of $N(x) \cup \{x\}$
of size $d$, and because the induced graph on $N(x) \cup \{x\}$ is
a $d$-complete hypergraph, this means that
$E'' = E' \cup \{z_{i-1},z_{i+1}\}$ is an edge of $\HH$.  The edge $E''$ is
distance 1 from $E_{i-1}$ and $E_{i+1}$.  We can replace $E$ in
the chain $E_0,\ldots,E_t$ with $E''$ and still have a proper chain of length
$t$ in $\HH\backslash E$ from $H$ to $H'$.  Moreover, this chain must be irredundant,
because if it was shorter, then this would give rise to shorter chain in $\HH$,
contradicting the fact that $t$ is the length of the shortest chain.
In case $(ii)$, let $z$ be the vertex of $E$ such that $\{z\} = E \setminus E_{i-1}$.
Then $E_{i-1} \setminus \{x\}$ and $z$ are in $N(x) \subseteq N(x) \cup \{x\}$.
Thus $E' = (E_{i-1} \setminus \{x\})\cup \{z\}$ is also an edge of $\HH$.  
Furthermore, $E'$ is distance one way from $E_{i-1}$ and
$E_{i+1}$  (because $z$ is added to $E$, we have $z \in E_{i+1}$).  So,
we can replace $E$ in the chain by $E'$ and get a chain of the correct
length in $\HH\backslash E$.  Finally, in case $(iii)$, let $z$ be the
vertex in $E$ such that $\{z\} = E \setminus E_{i+1}$.  Then 
$z \in N(x)$ and $(E_{i+1} \setminus\{x\}) \subseteq N(x)$.
This means that $E' = (E_{i+1}\setminus\{x\}) \cup \{z\}$
is an edge of $\HH$.  But this edge is distance one from both $E_{i-1}$
and $E_{i+1}$, so, as we did before, we can replace $E$ with $E'$ to 
get a chain of the desired length.

We can now show that $\HH \backslash E$ is also triangulated.
If the vertex $x \in E$ only appears in $E$, then $E$ is
a $v$-leaf.  Then $\HH \backslash E
= \HH \backslash \{x\} = \HH_{\X\backslash \{x\}}$, and it is clear that
$\HH_{\X \backslash \{x\}}$ is
a triangulated hypergraph.  So, suppose that
there are two or more edges that contain
$x$.  If $\Y \subseteq \X$
with $x \not\in \Y$, then the induced hypergraph of $\HH\backslash E$ on $\Y$
is the same as the induced hypergraph of $\HH$ on $\Y$, so there exists
a vertex $z \in \Y$ such that the induced hypergraph on $N(z) \cup \{z\}$ is a $d$-complete
hypergraph.  It remains to consider the case when $x \in \Y$.  Let $N_\Y(x)$ denote the neighbors
of $x$ in $(\HH\backslash E)_{\Y}$.  Note that $N_\Y(x)\cup \{x\} \subseteq N(x) \cup \{x\}$.  Since
the induced hypergraph on $N(x) \cup \{x\}$ is a $d$-complete hypergraph, any
induced subgraph on a subset of $N(x) \cup \{x\}$ is also a $d$-complete hypergraph.  So
the induced hypergraph $(\HH\backslash E)_{N_\Y(x) \cup \{x\}}$ is a $d$-complete hypergraph.
Thus $\HH \backslash E$ is triangulated.

Finally, by Lemma \ref{H'prop-con} we know that $\HH'$ is properly-connected.
The reason that $\HH'$ is triangulated follows from the fact that
\[ \HH' = \HH \backslash \{x,x_2,\ldots,x_d,z_1,\ldots,z_t\}
= \HH_{\X \backslash \{x,x_2,\ldots,x_d,z_1,\ldots,z_t\}}\]
where
$E = \{x,x_2,\ldots,x_d\}$ and $N(E) = \{z_1,\ldots,z_t\}$.
\end{proof}

We come to the main result of this section.

\begin{theorem}
Suppose that $\HH$ is a $d$-uniform triangulated hypergraph.  Then the
graded Betti numbers of $\I(\HH)$ can be computed recursively
using the formula
\[\beta_{i,j}(\I(\HH)) = \beta_{i,j}(\I(\HH\backslash E))
+ \sum_{l=0}^i \binom{t}{l}
\beta_{i-1-l,j-d-l}(\I(\HH'))\]
where $E$ is a splitting edge, $t = |N(E)|$, and
$\HH'$ and $\HH \backslash E$ are also $d$-uniform triangulated hypergraphs.
Here, $\beta_{-1,j}(\I(\HH')) = 1$ if $j =0$ and $0$ if $j \neq 0$.
\end{theorem}

\begin{proof} By Lemma \ref{triangulatedlemma}, the triangulated hypergraph
$\HH$ has a splitting edge $E$.  Furthermore, since both hypergraphs $\HH\backslash E$
and $\HH'$ are triangulated hypergraphs, they also have a splitting edges.  Thus,
by repeatedly using the formula of Theorem \ref{theoremsplitting} we get the
recursive formula.
\end{proof}

It is well known that the graded Betti numbers for an arbitrary monomial ideal
may depend 
upon the characteristic of $k$.  However, as a consequence of the above formula 
we obtain the following corollary.

\begin{corollary}\label{cor.hyperchar}
Suppose that $\HH$ is a triangulated hypergraph. Then the graded Betti numbers of $\I(\HH)$ are 
independent of the characteristic of the ground field and can be computed recursively.
\end{corollary}

When restricted to simple graphs, we get a particularly nice corollary.

\begin{corollary}\label{cor.graphchar}
Suppose that $G$ is a chordal graph.  Then the graded Betti numbers
of $\I(G)$ are independent of the characteristic of the ground field and can be computed recursively.
\end{corollary}

\noindent Jacques \cite{J} and Jacques and Katzman \cite{JK} first proved
Corollary \ref{cor.graphchar} in the special case that $G$ is a forest, a subclass of chordal graphs.


\section{Properly-connected hypergraphs and regularity}

In this section we investigate the Castelnuovo-Mumford regularity of
the edge ideal $\I(\HH)$ associated to a properly-connected hypergraph $\HH$.
For such a hypergraph, we bound $\reg(\I(\HH))$ below
by combinatorial invariants of the hypergraph.  When $\HH=G$ is a simple
graph, we also provided an upper bound.
 In the case
that $\HH$ is also triangulated, we explicitly compute $\reg(\I(\HH))$. Our exact formula 
for $\reg(\I(\HH))$
generalizes Zheng's formula \cite{Zheng2004} for the regularity of $\I(\HH)$
when $\HH = G$ is a forest.

We begin by relating the regularity of $\I(\HH)$ to
the regularity of edge ideals associated to sub-hypergraphs of $\HH$.  We produce
similar results for the projective dimension of $\I(\HH)$. We first make the convention 
that $\reg(0) = 1$ and if $\HH$ has no edges, we set $\pdim(\I(\HH)) = -1$.

\begin{lemma} \label{regL}
Let $E$ be any edge of a $d$-uniform properly-connected  hypergraph $\HH$ such that $\HH \backslash E$ 
is nonempty.
Let $t = |N(E)|$ and 
$\HH' = \{H \in \HH ~|~ \dist_{\HH}(H,E) \geq d+1\}$. If $L = (x^E) \cap \I(\HH\backslash E)$, then
\begin{enumerate}
\item[$(a)$] $\reg(L) = \reg(\I(\HH'))+d$, and
\item[$(b)$] $\pdim(L) = \pdim(\I(\HH')) + t.$
\end{enumerate}
\end{lemma}

\begin{proof}
We shall prove both results using Lemma \ref{betti intersection}.  For $(a)$
suppose $s = \reg(L)$.  So, there exists $a$ such that $\beta_{a,a+s}(L) \neq 0$.
By Lemma \ref{betti intersection}
\[\beta_{a+1-1,a+s}(L) = \sum_{l=0}^{a+1}\binom{t}{l}\beta_{a+1-1-l,a+s-d-l}(\I(\HH')).\]
Since every number in the summation on the right hand side is nonnegative, there
exists some $l$ such that $\beta_{a-l,a+s-d-l}(\I(\HH')) \neq 0$.  Hence,
$\reg(\I(\HH')) \geq s-d$, or equivalently, $\reg(\I(\HH')) + d \geq \reg(L)$.
Conversely, if $r = \reg(\I(\HH'))$, then there exists $b$ such that 
$\beta_{b,b+r}(\I(\HH'))\neq 0$.  But then since $b+r = (b+r+d)-d$, by
Lemma \ref{betti intersection}
\[0 \neq \beta_{b,(b+r+d)-d}(\I(\HH')) \leq
    \sum_{l=0}^{b+1}\binom{t}{l}\beta_{b+1-1-l,b+r+d-d-l}(\I(\HH')) = \beta_{b,b+r+d}(L).\]
So $\reg(\I(\HH'))+d \geq \reg(L) \geq \reg(\I(\HH'))+d$, as desired.

To prove $(b)$, suppose $N(E) = \{z_1,\ldots,z_t\}$.
In the proof of Lemma \ref{betti intersection} it was shown that
\[R/L \cong R_1/(z_1,\ldots,z_t) \otimes _k R_2/\I(\HH').\]
where $R_1 = k[z_1,\ldots,z_t]$ and $R_2 = k[x_1,\ldots,x_s]$, 
with $\{x_1,\ldots,x_s\} = \X \backslash N(E)$.  By tensoring
the resolutions of $R_1/(z_1,\ldots,z_t)$ and $R_2/\I(\HH')$, we get
\begin{eqnarray*}
\pdim(L)+1 = \pdim(R/L) &=& \pdim(R_1/(z_1,\ldots,z_t))+ \pdim(R_2/\I(\HH'))\\ 
&=& t + \pdim(R/\I(\HH')) = t + \pdim(\I(\HH'))+1.
\end{eqnarray*}
The desired identity is obtained by comparing the first and last values of the above equality.
\end{proof}

\begin{theorem}\label{reg pdim 2}
Let $E$ be any edge of a $d$-uniform properly-connected hypergraph $\HH$
such that $\HH \backslash E$ is nonempty.
Let $t =|N(E)|$.
Then
\begin{enumerate}
\item[$(a)$] $\reg(\I(\HH))\leq \max\{\reg(\I(\HH \backslash E)), \reg(\I(\HH')) + d -1\}.$
\item[$(b)$] $\pdim(\I(\HH))\leq \max\{\pdim(\I(\HH\backslash E)),\pdim(\I(\HH'))+t+1\}.$
\end{enumerate}
Furthermore, if $E$ is a splitting edge, then we have equality in both $(a)$ and $(b)$.
\end{theorem}

\begin{proof}
Set $L =  (x^E) \cap \I(\HH\backslash E)$.
The two inequalities then follow by using the short exact sequence
\[0 \rightarrow L \rightarrow (x^E) \oplus \I(\HH\backslash E)
\rightarrow \I(\HH) \rightarrow 0\]
and Lemma \ref{regL}
to bound $\reg(I(\HH))$ and $\pdim(\I(\HH))$, noting that since $\HH \backslash E$ is nonempty, 
$\reg(\HH \backslash E) \ge d$.
When $E$ is a splitting edge, the equalities are a result of
the formulas of Corollary \ref{reg pdim}.
\end{proof}

We now focus our attention on using combinatorial information from
$\HH$ to bound $\reg(\I(\HH))$.  More precisely, the regularity
will be expressed using the following terminology.

\begin{definition}
Let $\HH$ be a $d$-uniform properly-connected hypergraph.  Two edges $E,H$ of
$\HH$ are {\bf $t$-disjoint} if $\dist_{\HH}(E,H) \geq t$.  A set of edges $\E' \subseteq \E$
is {\bf pairwise $t$-disjoint} if every pair of edges of $\E'$ is $t$-disjoint.
(We thank Jeremy Martin
for suggesting this name.)
\end{definition}

\begin{remark}
When $\HH$ is a $d$-uniform properly-connected hypergraph, then two edges $E$ and $H$ are
$d$-disjoint if and only if $E \cap H = \emptyset$; that is, $E$ and $H$ are disjoint in the usual 
sense. When $\HH =G$ is a simple graph, Zheng's definition \cite[Definition 2.15]{Zheng2004} for two 
edges to be {\bf disconnected} is equivalent to our definition that the two edges be 3-disjoint in $G$.
\end{remark}

We come to the first main result of this section.

\begin{theorem} \label{regtheorem}
Let $\HH$ be a $d$-uniform properly-connected hypergraph. Then $\beta_{i-1,di}(\I(\HH))$ equals the 
number of sets of $i$ pairwise $(d+1)$-disjoint edges of $\HH$. In particular, if $c$ is the maximal 
number of pairwise $(d+1)$-disjoint edges of $\HH$ then 
$$\reg(\I(\HH)) \ge (d-1)c+1.$$
\end{theorem}

\begin{proof} The first statement of the theorem implies that $\beta_{c-1,dc}(\I(\HH)) \not= 0$.
 Thus, $dc-(c-1) \le \reg(\I(\HH))$ and the second statement is proved. We shall prove the first 
statement of the theorem. In the case $d=2$, this is the content of \cite[Lemma 2.2]{K}. We generalize 
Katzman's arguments to the more general situation.

Recall that $\E = \{E_1, \dots, E_s\}$ and let 
$\TT: 0 \rightarrow T_s \stackrel{\partial_s}{\rightarrow} \dots \stackrel{\partial_2}{\rightarrow} 
T_1 \stackrel{\partial_1}{\rightarrow} \I(\HH) \rightarrow 0$ be the Taylor resolution of $\I(\HH)$. 
Then $T_i$ is a free $R$-module with generators $e_{j_1, \dots, j_i}$, for 
$1 \le j_1 < \dots < j_i \le s$, and the boundary map $\partial_i$ is defined by
$$\partial_i(e_{j_1, \dots, j_i}) = 
\sum_{k=1}^i (-1)^k \mu_k e_{j_1, \dots, \widehat{j_k}, \dots, j_i},$$
where $\widehat{j_k}$ indicates the removal of $j_k$ and $\mu_k = x^{E_{j_k} 
\backslash (\cup_{l \not= k} E_{j_l})}$. Let $\m = (x_1, \dots, x_n)$ be the maximal homogeneous 
ideal in $R$. It is well-known that the graded Betti numbers of $\I(\HH)$ are given by 
$$\beta_{i-1,j}(\I(\HH)) = \dim_k H_i(\TT \otimes_R R/\m)_j.$$

Observe that generators of degree $di$ of $T_i$ are $e_{j_1, \dots, j_i}$'s where 
$E_{j_1}, \dots, E_{j_i}$ are pairwise disjoint. Consider one such generator 
$e_{j_1, \dots, j_i}$. Let $\HH_1$ be the induced sub-hypergraph of $\HH$ on the 
vertices in $\bigcup_{k=1}^i E_{j_k}$. It can be seen that for $1 \le k \le i$, $E_{j_k}$ is 
disjoint from $\bigcup_{l \not= k} E_{j_l}$ and hence, $\mu_k \in \m$. Thus, the image of 
$\partial_i(e_{j_1, \dots, j_i})$ in $\TT \otimes_R R/\m$ is 0. Also, if $\HH_1$ contains an 
edge $E_t$ different from $E_{j_1}, \dots, E_{j_i}$, then since 
$E_t \subseteq \bigcup_{k=1}^i E_{j_k}$, we have $\partial_{i+1}(e_{j_1, \dots, j_i, t}) = 
e_{j_1, \dots, j_i}$. That is, if $\HH_1$ contains an edge different from $E_{j_1}, \dots, E_{j_i}$ 
then the image of $e_{j_1, \dots, j_i}$ in $H_i(\TT \otimes_R R/\m)$ is 0. Furthermore, the image 
of $e_{j_1, \dots, j_i}$ in $\TT \otimes_R R/\m$ is in the image of $\partial_{i+1}$ if and only 
if it is the image of $\partial_{i+1}(e_{l_1, \dots, l_{i+1}})$, where 
$\{l_1, \dots, l_{i+1}\} = \{j_1, \dots, j_i\} \cup \{t\}$ for some $t$. This implies that in the 
expansion of $\partial_{i+1}(e_{l_1, \dots, l_{i+1}})$, we must have $\mu_t = 1$, i.e., 
$\HH_1$ contains the edge $E_t$ different from $E_{j_1}, \dots, E_{j_i}$. 

It remains to show that $E_{j_1}, \dots, E_{j_i}$ are pairwise disjoint edges of $\HH$ such that 
the induced sub-hypergraph of $\HH$ on the vertices of $\bigcup_{k=1}^i E_{j_k}$ contains no other 
edges if and only if $E_{j_1}, \dots, E_{j_i}$ are pairwise $(d+1)$-disjoint edges in $\HH$. 

Suppose first that $E_{j_1}, \dots, E_{j_i}$ are pairwise disjoint edges of $\HH$ such that the 
induced sub-hypergraph $\HH_1$ on the vertices of $\bigcup_{k=1}^i E_{j_k}$ contains no other edges. 
Clearly, since $E_{j_k} \cap E_{j_l} = \emptyset$ for $k \not= l$, we have 
$\dist_{\HH}(E_{j_k},E_{j_l}) \ge d$. Now, suppose there exist $k \not= l$ so that 
$\dist_{\HH}(E_{j_k},E_{j_l}) = d$. Then there is a proper chain $E_{j_k} = F_0, F_1, \dots, F_d = 
E_{j_l}$. By Lemma \ref{lem: chains}, the vertices of $F_1$ are in $E_{j_k} \cup E_{j_l}$, so $F_1$ is 
an edge in $\HH_1$. This implies that $F_1$ has to be one of the $\{E_{j_1}, \dots, E_{j_i}\} 
\backslash \{E_{j_k}, E_{j_l}\}$. This is a contradiction since $F_1 \cap E_{j_k} \not= \emptyset$.

Conversely, suppose that $E_{j_1}, \dots, E_{j_i}$ are pairwise $(d+1)$-disjoint edges of $\HH$. 
Let $\HH_1$ be the induced sub-hypergraph of $\HH$ on the vertices of $\bigcup_{k=1}^i E_{j_k}$. 
By contradiction, suppose $\HH_1$ contains an edge $E$ different from $E_{j_1}, \dots, E_{j_i}$. 
Then $E \subseteq \bigcup_{k=1}^i E_{j_k}$. Without loss of generality, we may assume that $E \cap 
E_{j_1} \not= \emptyset$. Then there is a proper chain $E_{j_1} = F_0, F_1, \dots, F_l = E$ for 
some $l < d$. By Lemma \ref{lem: chains}, the vertices of $F_1$ are in $E_{j_1} \cup E$. Thus, 
$F_1$ is also an edge of $\HH_1$. This implies that there exists $j_k \not= j_1$ so that $F_1$ has a 
nonempty intersection with $E_{j_k}$ (otherwise, $F_1 \subseteq E_{j_1}$, which is a contradiction). 
However, we now have $\dist_{\HH}(F_1, E_{j_k}) = d - |F_1 \cap E_{j_k}| \le d-1$, whence 
$\dist_{\HH}(E_{j_1},E_{j_k}) \le d$, which is again a contradiction.
\end{proof}

When $\HH$ is a graph we also obtain an especially appealing upper bound for the regularity of $\I(\HH)$.

\begin{definition}
Let $\HH = (\X,\E)$ be a hypergraph.  A {\bf matching} of $\HH$ is defined to be a
subset $\E' \subseteq \E$ consisting of pairwise disjoint edges.
The {\bf matching number} of $\HH$, denoted $\alpha'(\HH)$, is the largest size
of a maximal matching in $\HH$.
\end{definition}


\begin{theorem} \label{cor.matching}
Let $G$ be a finite simple graph.  Then
$$\reg(R/\I(G))\leq \alpha'(G)$$
where $\alpha'(G)$ is the matching number of $G$.
\end{theorem}

\begin{proof} It can be seen from the Taylor resolution that
$$\reg(\I(G)) \le 
\max \{ \deg \lcm(x^{E_1}, \dots, x^{E_i}) - i ~|~ \{E_1, \dots, E_i\} \subseteq \E \} + 1.$$
Since any edge of $G$ has 2 vertices, it can be seen that 
$i+1 \le \deg \lcm(x^{E_1}, \dots, x^{E_i}) \le 2i$. Let $\deg \lcm(x^{E_1}, \dots, x^{E_i}) = i+k$ 
for some $1 \le k \le i$. It suffices to show that we can always find a matching of size 
$k$ among $\{E_1, \dots, E_i\}$. To this end, we shall use induction on $i+k$.

If $i+k = 2$, i.e., $i = k = 1$, then the statement is clear. Suppose now that 
$i+k > 2$. If $k = 1$ or $k=i$ then the statement is also clear. 
Assume that $1 < k < i$. If $E_i$ is disjoint from $E_j$ for all $j < i$, 
then $\deg \lcm(x^{E_1}, \dots, x^{E_{i-1}}) = i+k-2 = (i-1) + (k-1)$. By induction, there exists a 
matching $S \subset \{E_1, \dots, E_{i-1}\}$ of size $(k-1)$. It is easy to see that $S \cup \{E_i\}$ 
is now a matching of size $k$. It remains to consider the case that at least a vertex of $E_i$ is 
also a vertex of $E_j$ for some $j < i$. In this case, we have 
$\deg \lcm(x^{E_1}, \dots, x^{E_{i-1}}) \ge i+k-1 = (i-1)+k$. By induction, there is a 
matching $S \subset \{E_1, \dots, E_{i-1}\}$ of size $k$, and the statement is proved.
\end{proof}

Theorem \ref{cor.matching} seems to give an interesting bound for the regularity of edge ideals 
with a simple proof which may have been overlooked.

When $\HH$ is a triangulated hypergraph, the lower bound of Theorem \ref{regtheorem} turns out to be the exact formula for the regularity of $\I(\HH)$.

\begin{theorem}\label{regularitytheorem}
Suppose that $\HH$ is a $d$-uniform properly-connected triangulated hypergraph,
If $c$ is the maximum number of pairwise $(d+1)$-disjoint edges of $\HH$, then
\[\reg(\I(\HH)) = (d-1)c + 1.\]
\end{theorem}

\begin{proof}  The proof is similar to the one given by \cite{HaVanTuyl2005} in the case for forests.  
We proceed by  induction on the number of edges of $\HH$.
If $\HH$ only has one edge $E$, then $\I(\HH) = (x^E)$.
Because $\I(\HH)$ is principal, it is clear that $\reg(\I(\HH)) = d$.
But then it is clear that the formula holds since $1$ is the maximal number
of pairwise $(d+1)$-disjoint edges.

So, suppose $\HH$ has at least two edges. Since $\HH$
is triangulated, by Lemma \ref{triangulatedlemma} there is a splitting edge $E \in \HH$ ($\HH \backslash E$ is nonempty in this case) such that $\HH\backslash E$ and $\HH'$ are also
$d$-uniform properly-connected triangulated hypergraphs.
Since $E$ is a splitting edge, by Corollary \ref{reg pdim 2} we have
\[\reg(\I(\HH)) = \max\{\reg(\I(\HH\backslash E)), \reg(\I(\HH'))+d-1\}.\]
By induction $\reg(\I(\HH\backslash E)) = (d-1)c_1+1$ where $c_1$ is the maximal number
of pairwise $(d+1)$-disjoint edges of  $\HH\backslash E$, and
$\reg(\I(\HH')) = (d-1)c_2+1$
where $c_2$ is the maximal number of pairwise $(d+1)$-disjoint edges of $\HH'$.
So
\[\reg(\I(\HH)) = \max\{(d-1)c_1+1,(d-1)c_2+d\}.\]
If we let $c$ denote the maximal number of pairwise $(d+1)$-disjoint edges of $\HH$,
then since $(d-1)c_2 + d = (d-1)(c_2 + 1) + 1$
to complete the proof it suffices for us to show that $c = \max\{c_1,c_2+1\}$.

Let $\mathcal{E}_1$ be the set of the $c_1$ pairwise $(d+1)$-disjoint edges of $\HH\backslash E$. The
edges of $\mathcal{E}_1$ are also a set of pairwise $d+1$-disjoint edges of $\HH$.
To see this fact, suppose that $H,H'$ are two $(d+1)$-disjoint edges in $\HH\backslash E$
that are not $(d+1)$-disjoint in $\HH$.  That is, $\dist_{\HH}(H,H') \leq d$.  But
because $H \cap H' = \emptyset$, we must have $\dist_{\HH}(H,H') = d$.
Let $H = E_0,\ldots,E_d=H'$ be the proper irredundant chain of length $d$ in $\HH$.
Since this chain is not in $\HH\backslash E$, we must have $E = E_i$ for some $i=\{1,\ldots,d-1\}$.
Consider the edges $E_{i-1}$ and $E_{i+1}$ in the chain that
occur before and after, respectively, the edge $E$.  The splitting edge $E$ of Lemma
\ref{triangulatedlemma} is picked so that it contains a vertex $x$ such
that the induced graph on $N(x) \cup \{x\}$ is a $d$-complete hypergraph.
We can now adapt the proof given in Lemma \ref{triangulatedlemma}
that showed that $\HH\backslash E$ was properly-connected to show that $E$ can be replaced
by an edge $E' \in \HH\backslash E$.  As a consequence,  we get a path of length $d$ from $H$ to $H'$
in $\HH\backslash E$.
But this contradicts the fact that $\dist_{\HH\backslash E}(H,H') \geq d+1$.
Thus $|\mathcal{E}_1| = c_1 \leq c$.  

If $\mathcal{E}_2$
is a set of $c_2$ pairwise $(d+1)$-disjoint edges of $\HH'$, we claim that
$\mathcal{E}_2 \cup \{E\}$ is a set of pairwise $(d+1)$-disjoint edges of $\HH$.
Indeed, for any edge $H \in \HH'$, $\dist_{\HH}(E,H) > d$,
and so in particular, $E$ and $H$ is $(d+1)$-disjoint for every edge $H \in \mathcal{E}_2$.
Thus $|\mathcal{E}_2 \cup \{E\}| = c_2 + 1 \leq c$. Thus $c \geq \max\{c_1,c_2+1\}$.

Suppose that  $c > \max\{c_1,c_2+1\}$.  Let $\mathcal{E}$ be a set of $c$
pairwise $(d+1)$-disjoint edges of $\HH$.  If $E \not\in\mathcal{E}$,
then $\mathcal{E}$ is also a set of pairwise $(d+1)$-disjoint edges of $\HH\backslash E$,
and so $c = |\mathcal{E}| \leq c_1$, a contradiction.  If $E \in \mathcal{E}$,
then $\mathcal{E} \backslash \{E\}$ is a set of  pairwise $(d+1)$-disjoint
edges of $\HH'$ since any other edge $H \in \mathcal{E}$
must have $\dist_{\HH}(E,H) > d$.  But this would imply that $c-1 \leq c_2$,
again a contradiction.  Hence $c = \max\{c_1,c_2+1\}$.
\end{proof}

Theorem \ref{regularitytheorem} gives the following interesting corollary for simple graphs, which was first proved by Zheng \cite{Zheng2004} in the special case that $G$ was a forest.

\begin{corollary} \label{cor.regchordal}
Suppose that $G$ is chordal graph.  If $c$ is the maximum number of
pairwise $3$-disjoint edges of $G$, then
\[\reg(\I(G)) = c + 1.\]
\end{corollary}

\begin{example}  \label{upperbound}
The bounds for the regularity in Theorems \ref{regtheorem} and Theorem \ref{cor.matching}
are sharp.  If $\HH$ is any triangulated hypergraph, then the lower bound in Theorem \ref{regtheorem}
is achieved by Theorem \ref{regularitytheorem}.  To show that
the upper bound in Theorem \ref{cor.matching} is achieved, consider the the edge ideal of $C_5$, the
five-cycle.  So $\I(G) = (x_1x_2,x_2x_3,x_3x_4,x_4x_5,x_5x_1)$.  Then
$\alpha'(G) = 2$ (for example, take edges $E_1 = x_1x_2$ and $E_2 = x_3x_4$).  So $\reg(\I(G)) \leq 3$.  
In fact we have equality since the resolution of $\I(G)$ is
\[0 \rightarrow R(-5) \rightarrow R^5(-3) \rightarrow R^5(-2) \rightarrow \I(G) \rightarrow 0. \]
\end{example}

In the study of squarefree monomial ideals, the theory of Alexander duality has proved to be 
significant in many ways. We round out this section by
relating some algebraic invariants of 
edge ideals and their Alexander duals.

\begin{definition} \label{Alexander dual}
Let $I = (x_{11} \cdots x_{1i_1}, \dots, x_{r1} \cdots x_{ri_r}) \subseteq k[x_1, \dots, x_n]$ be 
a squarefree monomial ideal. Then the {\bf Alexander dual} of $I$ is defined to be
$$I^\vee = (x_{11}, \dots, x_{1i_1}) \cap \dots \cap (x_{r1}, \dots, x_{ri_r}).$$
\end{definition}

\begin{definition} \label{vertex cover}
Let $G$ be a graph. A subset $V$ of the vertices of $G$ is called a {\bf vertex cover} if every 
edge in $G$ is incident to at least a vertex in $V$; a {\bf minimal vertex cover} is a vertex cover
$V$ with the property that no proper subset of $V$ is. 
The smallest size of a minimal vertex cover of $G$ is denoted by $\nu(G)$. 
The graph $G$ is {\bf unmixed} if all its minimal vertex covers have the same cardinality $\nu(G)$.
\end{definition}

\begin{remark} \label{Alexander gens}
The operation of taking the Alexander dual of a squarefree monomial ideal brings generators to 
primary components. The minimal generators of $\I(G)^\vee$ correspond to 
minimal vertex covers of $G$.
\end{remark}

\begin{theorem} \label{thm.Konig}
Let $G$ be a simple graph.
\begin{enumerate}
\item If $G$ is unmixed, then
\[\reg(\I(G)) \le \height \I(G) + 1 \le \reg(\I(G)^\vee)+1
~\text{and}~ \pdim(\I(G)^{\vee}) \leq \height \I(G) \leq \pdim(\I(G))+1.\]
\item If $G$ is not unmixed, then 
\[\reg(\I(G)) \le \height \I(G) +1\le \reg(\I(G)^\vee)
~\text{and}~\pdim(\I(G)^{\vee}) \le \height \I(G) \le \pdim(\I(G)).\]
\end{enumerate}
\end{theorem}

\begin{proof} It suffices to prove the inequalities involving the regularity,
since the bounds on the projective dimension follow from the
identities  $\reg(\I(G)) = \pdim(R/\I(G)^\vee)$ and 
$\reg(\I(G)^\vee) = \pdim(R/\I(G))$ (see, for example, \cite[Theorem 5.59]{MillerSturmfels2004}).
Observe that if $\E'$ is a matching in $G$ then any vertex cover must contain at least a vertex of every edge in $\E'$. Thus, $\alpha'(G) \le \nu(G) = \height \I(G)$. It follows from Theorem \ref{cor.matching} that $\reg(\I(G)) \le \height \I(G) + 1$. Since $\nu(G)$ is the least 
generating degree of $\I(G)^\vee$, we have $\nu(G) \le \reg(\I(G)^\vee)$ and thus (1) follows. To 
prove (2) observe that when $G$ is not unmixed, $\reg(\I(G)^\vee)$ is at least the largest 
generating degree of $\I(G)^\vee$, which is at least $\nu(G)+1$.  
\end{proof}


\section{Properly-connected hypergraphs and linear first syzygies}

In \cite{Fr} Fr\"oberg gave a characterization of edge ideals
of simple graphs with linear resolutions.  In this section,
we obtain a partial generalization of Fr\"oberg's result to the
class of properly-connected hypergraphs.  Specifically,
we describe when $\I(\HH)$ has linear first syzygies.

Let us first recall Fr\"oberg's result.  If $G$ is a simple graph,
then the {\bf complement of $G$}, denoted $G^c$, is the
graph whose vertex set is the same as $G$, but
whose edge set is defined by the rule $E \in G^c$ if and only
$E \not\in G$.  Fr\"oberg then showed:

\begin{theorem} \label{froberg}
Let $G$ be a simple graph.  Then
$\I(G)$ has a linear resolution if and only if $G^c$ is a chordal graph.
\end{theorem}

When $\HH$ is a $d$-uniform properly-connected hypergraph,
we define the {\bf complement of $\HH$}, denoted $\HH^c$, as
\[\HH^c = \{E \subseteq \X ~\big|~ |E| = d ~~\text{and}~~ E \not\in \HH\}.\]
So, one might expect Theorem \ref{froberg}
generalizes to $d$-uniform properly-connected hypergraphs as follows:
$\I(\HH)$ has a linear resolution if and only if $\HH^c$ is a triangulated
hypergraph.  Unfortunately, this is not the case, as shown below,
since $\HH^c$
need not be properly-connected.

\begin{example}
Let $\X = \{x_1,x_2,x_3,x_4,x_5\}$.  Let $\HH = \K^3_5 \backslash \{x_1x_2x_3,x_3x_4x_5\}$,
i.e., $\HH$ is the $3$-uniform complete hypergraph of order $5$ with
two edges removed.  Then $\HH^c = \{x_1x_2x_3,x_3x_4x_5\}$
is not properly-connected since the two edges intersect at $x_3$, but
there is no properly-irredundant chain of length 2 between the two edges.
Because $\HH^c$ is not even properly-connected, the notion of a triangulated
hypergraph is undefined.  However, the ideal $\I(\HH)$ has the linear
resolution
\[0 \rightarrow R^4(-5) \rightarrow R^{11}(-4) \rightarrow  R^8(-3) \rightarrow \I(\HH) \rightarrow
0.\]
\end{example}

We take the first step towards generalizing Theorem \ref{froberg}
by asking when $\I(\HH)$ must have linear first syzygies.  Like
our previous results, the
distance between edges is key.

\begin{definition}
The {\bf edge diameter} of a $d$-uniform properly-connected
hypergraph $\HH$ is
\[ \operatorname{diam}(\HH) = \max \{ \dist_{\HH}(E,H) \mid E, H \in \HH\},\]
where the diameter is infinite if there exist two edges not connected by any proper chain.
\end{definition}

Since $\I(\HH)$ is a monomial ideal, we know that its first syzygy module is generated by
syzygies $S(x^E, x^H)$, for $E, H \in \E$.
Moreover, it is clear that $S(x^E, x^H)$ is a linear syzygy if and only if
$\dist_{\HH}(E,H) = 1$.  We shall see that these syzygies generate all of the syzygies on $\I(\HH)$ if
the diameter of $\HH$ is small enough.  Indeed, a short enough proper chain will give us a way of
writing $S(x^E, x^H)$ as a telescoping sum of linear syzygies.  The next theorem
generalizes \cite[Theorem 3.17]{Zheng2004}.  

\begin{theorem}\label{linearsyzygies}
Suppose that $\HH$ is a $d$-uniform properly-connected hypergraph.
Then $\I(\HH)$ has linear first syzygies if and only if diam$(\HH) \leq d$.
\end{theorem}

\begin{proof}
Assume first that diam$(\HH) \le d$. It follows from the Taylor resolution that the first syzygy
module of $\I(\HH)$ is generated by syzygies $S(x^E,x^H)$, where $E,H \in \mathcal{E}$.
We shall show that $S(x^E,x^H)$ is generated by linear syzygies. Let $t = \dist_\HH(E,H)$. Then, 
since diam$(\HH) \le d$, we have $t \le d$.
If $(E_0,\ldots,E_t)$ is the proper irredundant chain, then
by Lemma \ref{lem: chains} we can write
$E = E_0= \{ z_1, \ldots, z_d\},$ $E_i = \{y_1, \ldots, y_i, z_{i+1}, \ldots, z_d\}$ where
$y_i \notin E_j$ for $j < i,$ and $E_t = H$.

It can be seen that $S(x^E,x^H)$ is given by the
equality $y_1 \cdots y_t x^{E_0} - z_1\cdots z_t x^{E_t} = 0.$
Furthermore,
$$y_1 \dots y_tx^{E_0} - z_1 \dots z_t x^{E_t} =
\sum_{k=0}^{t-1} \left(\prod_{i=1}^k z_i \prod_{j=k+2}^t y_j\right)
(y_{k+1}x^{E_k} - z_{k+1}x^{E_{k+1}}).$$
Thus, $S(x^E,x^H)$ is generated by linear syzygies.

Conversely, suppose that $\I(\HH)$ has linear first syzygies, that is,
$\beta_{1,j}(\I(\HH)) = 0$ for $j \neq d+1$.  If
$\operatorname{diam}(\HH) \geq d+1$, then this implies that there exists at least
two edges $E,H$ with $\dist_{\HH}(E,H) \geq d+1$, i.e., $\{E,H\}$ is
a set of pairwise $(d+1)$-disjoint edges of $\HH$.  By Theorem \ref{regtheorem}
this implies that $\beta_{1,2d}(\I(\HH)) \neq 0$.  But this contradicts
the fact that $\I(\HH)$ has linear first syzygies.
\end{proof}

\begin{example} If $\operatorname{diam}(\HH) \leq d$ is small, 
$\I(\HH)$ may still have nonlinear second syzygies.  For example,
if $G = C_5$ is the 5-cycle, then $\operatorname{diam}(G) = 2$.  However
$\I(G) = (x_1x_2,x_2x_3,x_3x_4,x_4x_5,x_5x_1)$ has nonlinear
second syzygies since $\beta_{2,5}(\I(G)) = 1$, as shown in Example
\ref{upperbound}.
\end{example}

Interestingly, if $\HH$ is triangulated, knowing that $\I(\HH)$ has linear
first syzygies is enough to know that the entire resolution
of $\I(\HH)$ is linear.

\begin{corollary} \label{linsyzcor}
Suppose that $\HH$ is a $d$-uniform properly-connected hypergraph that is
also triangulated.  Then the following are equivalent:
\begin{enumerate}
\item[$(a)$] $\I(\HH)$ has a linear resolution.
\item[$(b)$] $\I(\HH)$ has linear first syzygies.
\item[$(c)$] $\operatorname{diam}(\HH) \leq d$.
\end{enumerate}
\end{corollary}

\begin{proof}
The implication $(a) \Rightarrow (b)$ is immediate, and $(b) \Rightarrow (c)$
is a consequence of Theorem \ref{linearsyzygies}.  To show that $(c)
\Rightarrow (a)$,  the bound on $\operatorname{diam}(\HH)$ implies that $\HH$
cannot have two or more
pairwise $(d+1)$-disjoint edges (otherwise diam$(\HH) > d$).  By
Theorem \ref{regularitytheorem} this implies that $\reg(\I(\HH)) = (d-1)+1 =d$.
Since $\I(\HH)$ is generated in degree $d$, this forces $\I(\HH)$ to have a linear
resolution.
\end{proof}

Restricted to simple graphs, Corollary \ref{linsyzcor} gives the following result.

\begin{corollary} Suppose that $G$ is a chordal graph. Then the following are equivalent:
\begin{enumerate}
 \item[$(a)$] $\I(G)$ has a linear resolution.
 \item[$(b)$] $\I(G)$ has linear first syzygies.
 \item[$(c)$] $\operatorname{diam}(G) \leq 2$.
\end{enumerate}
\end{corollary}


\section*{ Acknowledgements} The authors would especially like to thank Jessica Sidman
who made some contributions to this paper in its preliminary stages, and with whom
we had many useful discussions. The authors would like to thank J. Herzog and X. Zheng
for stimulating discussions on the regularity of edge ideals. Part of this research was carried out while
the second author visited the first at Tulane University. The authors acknowledge the support from 
Louisiana Board of Regents for this visit. The second author also thanks Tulane University 
for its hospitality during his visit. The second author further acknowledges the research 
support received by
NSERC. The computer algebra system {\tt CoCoA} \cite{Co} was used to generate examples.
We would also like to thank the two anonymous referees for their suggestions and comments. 


\end{document}